\documentstyle[12pt]{article}
\newcommand{\gc}[2]{\left[{#1}\atop{#2}\right]_q}
\newcommand{\qp}[4]{\left[{#1}\atop{#2},{#3},{#4}\right]_q}
\newcommand{\af}[5]{\mathop{{#1}(^{#2}_{{#3},{#4},{#5}})}}
\newcommand{\st}[2]{\left[{#1}\atop{#2}\right]}
\newcommand{\stt}[2]{\left\{{#1}\atop{#2}\right\}}
\newcommand{\eu}[2]{\left\langle{#1}\atop{#2}\right\rangle}

\newcommand{\be}{\begin{equation}}
\newcommand{\ee}{\end{equation}}
\begin{document}

\centerline{\Large\bf A probabilistic approach to $q$-polynomial}
\centerline{\Large\bf coefficients, Euler and Stirling numbers}
\bigskip

\centerline{\bf Alexander I. Il'inskii}
\bigskip

\centerline{\small Kharkov National University, Department of
Mathematics,} \centerline{\small 4 Svobody Sqr., 61077, Kharkov,
Ukraine} \centerline{\it e-mail:
Alexander.I.Iljinskii@univer.kharkov.ua}
\bigskip
\centerline{\bf Abstract}

It is known that Bernoulli scheme of independent trials with two
outcomes is connected with the binomial coefficients. The aim of
this paper is to indicate stochastic processes which are connected
with the $q$-polynomial coefficients (in particular, with the
$q$-binomial coefficients, or the Gaussian polynomials), Stirling
numbers of the first and the second kind, and Euler numbers in a
natural way. A probabilistic approach allows us to give very
simple proofs of some identities for these coefficients.
\medskip

\medskip

\bigskip
\section{Introduction} \setcounter{equation}{0}

The connection of binomial coefficients $\displaystyle\left(n\atop
k\right)$ with Bernoulli scheme of independent trials is well
known. To be more specific, let $x\in(0,1)$, $n$ be a positive
integer, ${\cal E}$ be a trial with two outcomes $0$ and $1$ with
probabilities in a single trial being equal to $1-x$ and $x,$
respectively. Let ${\cal E}$ be repeated $n$ times under the
condition that every outcome of any trial is independent of
outcomes of all other trials. Let ${\bf 1}(^n_k)$ denote an event
such that the outcome $1$ has happened $k$ times in $n$
repetitions of the trial ${\cal E}$. Then the probability $P({\bf
1}(^n_k))$ of this event equals
$$
P({\bf 1}(^n_k))={n\choose k}x^k(1-x)^{n-k}\,.
$$

The aim of this paper is to construct random processes which are
connected with $q$-polynomial coefficients \be \label{1.1}
\left[{i_1+i_2+\ldots +i_m}\atop{i_1, i_2,\ldots
,i_m}\right]_q:={{(q)_{i_1+i_2+\ldots
+i_m}}\over{(q)_{i_1}(q)_{i_2}\ldots (q)_{i_m}}}\,, \ee (where
$(q)_j:=(1-q)(1-q^2)\ldots(1-q^j)$ for $j\in{\bf N}$, $(q)_0:=1$),
in particular, $q$-binomial coefficients (or Gaussian
polynomials). The probabilistic approach gives very simple proofs
of some identities for $q$-polynomial coefficients.

We apply an analogues approach to  Stirling numbers of the first
and the second kind $\displaystyle\st{n}{k}$,
$\displaystyle\stt{n}{k}$ and Euler numbers
$\displaystyle\eu{n}{k}$ (we follow notation used in [3], Chapter
6).

In this paper, we introduce numbers associated with Euler and
Stirling numbers and obtain some identities for them.

The paper is organized as follows. In Section~2, we consider
$q$-binomial coefficients. We consider this special case of
$q$-polynomial coefficients separately for the convenience of
readers.

Section~3 is devoted to $q$-polynomial coefficients. We consider
the case $m=3$ in (\ref{1.1}) only. The general case is obvious
after that.

In Sections~4,~5,~6, we consider Stirling numbers of the second
and the first kind and Euler numbers, respectively.

In the beginning of each section we define a probability space or
a space equipped with a weight which is connected with
corresponding coefficients in a natural way (see Theorems~2.1,
3.1, 4.3, 5.2, 6.2). Using probabilistic arguments we give very
simple proofs of some identities for these coefficients (see
Theorems~2.2, 2.3, 3.2, 3.3, 4.4, 5.3, 6.3). We also introduce a
notion of coefficients associated with Stirling and Euler numbers
and deduce some identities for them (see Theorems~4.5, 5.4, 6.4).

This study has been stimulated by the paper~[4].

As usual, ${\bf N}=\{1, 2, 3,\ldots\}$ and ${\bf N_0}=\{0, 1, 2,
3,\ldots\}$  denote the sets of positive integers and nonnegative
integers, respectively. We use the following notation (see [1],
Chapter 3):
$$
\begin{array}{rcl}
(x;q)_n&:=&(1-x)(1-qx)\cdots (1-q^{n-1}x)\,\,\,\mbox{ for }
n\in{\bf
N}\,,\quad (x;q)_0:=1\,;\\
(q)_n&:=&(q;q)_n\,.
\end{array}
$$
The identity \be \label{1.2} (x;q)_m\cdot (q^mx;q)_n=(x;q)_{m+n}
\ee will be useful in the sequel.

We recall  that Stirling numbers of the first and the second kind
$\displaystyle\st{n}{k}$, $\displaystyle\stt{n}{k}$, and Euler
numbers $\displaystyle\eu{n}{k}$ may be defined for $n\in{\bf
N_0}$ and integer $k$ such that $0\leq k\leq n$ as numbers which
equal $1$ if $n=k=0$, and $0$, if $k<0$ or $k>n$, and satisfy the
following recurrence identities (see [3], Section~6.1)
\begin{eqnarray}
\label{1.3}
\st{n}{k}&=&\st{n-1}{k-1}+(n-1)\st{n-1}{k}\,,\nonumber\\
\stt{n}{k}&=&\stt{n-1}{k-1}+k\stt{n-1}{k}\,,\\
\eu{n}{k}&=&(n-k)\eu{n-1}{k-1}+(k+1)\eu{n-1}{k}\,,\nonumber
\end{eqnarray}
respectively.

We remind some elementary probabilistic concepts which we use in
this paper (see, for example, [2]). A {\it probability space} is a
pair $(\Omega, p)$ where $\Omega$ is a finite set (it is called a
{\it sample space}\,) and $p$ is a function on $\Omega$ (it is
called a {\it probability}\,) such that $p(\omega)\geq 0$ for all
$\omega\in\Omega$ and $\sum_{\omega\in\Omega}p(\omega)=1$. Every
set $A\subset\Omega$ is said to be an {\it event}\,. The {\it
probability of an event} $A$ is defined by $P(A)=\sum_{\omega\in
A}p(\omega)$. If $A, B\subset\Omega$ and $P(B)\neq 0$, then
$P(A|B):=P(A\cap B)/P(B)$ is said to be the {\it conditional
probability} of $A$ given $B$. The formula $P(A\cap B)=P(B)P(A|B)$
is called {\it the multiplication theorem of probability}\,. We
say that events $A_l$, $l=1,2,\ldots L$, form a {\it partition} of
a sample space $\Omega$ if $A_l$ are mutually exclusive and
$\cup_{l=1}^LA_l=\Omega$. In that case, {\it the formula of total
probability} is valid:
$$
P(B)=\sum\limits_{l=1}^LP(B|\,A_l)P(A_l)\,.
$$

We  also use the following notation:
$i_{(k)}:=\underbrace{i,i,\ldots,i}_{k}$ for $k\in{\bf N}$.

\bigskip

\section{$q$-binomial coefficients (Gaussian polynomials)}
\setcounter{equation}{0}
\bigskip

Let $x$ and $q$ be arbitrary real numbers from the interval
$(0,1)$, $n$ be a positive integer.

We will construct a special stochastic process. Let ${\cal E}$ be
a trial with two outcomes $0$ and $1$. We consider a sequence of
$n$ trials ${\cal E}$. We take the probability of outcomes $0$ and
$1$ in the first trial $\cal E$ to be equal $1-x$ and $x$,
respectively. Assume that the trial $\cal E$ is repeated $m$ times
and the outcome $0$ has occurred $j$ times where $0\leq j\leq m$.
Then we take the probability of $0$ and $1$ in the $(m+1)^{th}$
repetition of $\cal E$ to be $1-q^jx$ and $q^jx$, respectively.

A mathematical model of this sequence of $n$ trials ${\cal E}$ is
the probability space $(\Omega_n, p_{x,q,n})$ where
$$
\Omega_n:=\{\omega=(\varepsilon_1, \varepsilon_2,\ldots ,
\varepsilon_n): \varepsilon_k=0 \mbox{ or } 1, k=1, 2,\dots ,
n\}\,,
$$
the probability $p_{x,q,n}(\omega)$ of an elementary event
$\omega=(\varepsilon_1, \varepsilon_2,\ldots , \varepsilon_n)$ is
equal to \be \label{2.1} p_{x,q,n}(\omega)=
p_{x,q,n}((\varepsilon_1, \varepsilon_2,\ldots , \varepsilon_n))=
f_1\cdot f_2\cdot\ldots\cdot f_n\,, \ee where $f_1=x$, if
$\varepsilon_1=1$, $f_1=1-x$, if $\varepsilon_1=0$, and for every
integer $m$ ($1\leq m\leq n-1$), \be \label{2.50}
f_{m+1}=\cases{q^jx\,,& if $\varepsilon_{m+1}=1$\,,\cr 1-q^jx\,, &
if $\varepsilon_{m+1}=0$\,,\cr} \ee where $j=\#\{l:1\leq l\leq
m\,,\varepsilon_l=0\}$. The probability of an event
$A\subset\Omega_n$ is defined by \be \label{2.2}
P_{x,q,n}(A):=\sum\limits_{\omega\in A}p_{x,q,n}(\omega)\,. \ee
For the sake of brevity, sometimes we  write $P$ instead of
$P_{x,q,n}$. It is not difficult to see that
$P_{x,q,n}(A)=P_{x,q,m}(A)$ if $m<n$ and the event
$A\subset\Omega_n$ does not depend on the repetitions of the trial
${\cal E}$ with numbers $m+1, m+2,\ldots ,n$, that is $A$
satisfies the following condition: if
$(\varepsilon_1,\ldots,\varepsilon_m,\varepsilon_{m+1},\ldots,\varepsilon_n)\in
A$, then
$(\varepsilon_1,\ldots,\varepsilon_m,\delta_{m+1},\ldots,\delta_n)\in
A$ for all $\delta_{m+1},\ldots,\delta_n=0 \mbox{ or }1$.

Let us introduce some notation. For $i=0, 1$ and integers $m$
($1\leq m\leq n$) and $k$ ($0\leq k\leq n$) we define \be
\label{2.3} {\bf
i}^{(m)}:=\{\omega=(\varepsilon_1,\varepsilon_2,\ldots
\varepsilon_n)\in\Omega_n:\varepsilon_m=i\}\,, \ee \be \label{2.4}
{\bf i}(^n_k):=\{\omega=(\varepsilon_1,\varepsilon_2,\ldots
\varepsilon_n)\in\Omega_n:\#\{l:1\leq l\leq n,
\varepsilon_l=i\}=k\}\,. \ee

Using this notation we may write:
$$
P_{x,q,n}({\bf 0}^{(1)})=1-x\,, \qquad P_{x,q,n}({\bf 1}^{(1)})=x
$$
for every $n\geq 1$;
$$
P_{x,q,m+1}({\bf 1}^{(m+1)}|\,{\bf 1}(^{m}_k))=q^{m-k}x\,, \quad
P_{x,q,m+1}({\bf 0}^{(m+1)}|\,{\bf 1}(^{m}_k))=1-q^{m-k}x
$$
for every $m\geq 1$. It should be pointed out that the exponent of
the power $q^{m-k}$ is equal to the number of all trials which
lead to outcome $0$.

We give now the probabilities of some elementary events $\omega$
such that $\omega\in{\bf 1}(^n_k)$:
\begin{eqnarray*}
p((1_{(k)},0_{(n-k)}))&=&
x^k(x;q)_{n-k}\,,\\
p((0_{(n-k)},1_{(k)}))&=&
x^k(x;q)_{n-k}q^{k(n-k)}\,,\\
p((0,1,0,1_{(k-1)}, 0_{(n-k-2)}))&=& x^k(x;q)_{n-k}q^{2k-1}\,.
\end{eqnarray*}
It is not difficult to see that $p(\omega)=x^k(x;q)_{n-k}q^r$ for
every $\omega\in{\bf 1}(^n_k)$, where $r$ is an integer such that
$0\leq r\leq k(n-k)$. Therefore, we may write for all $n\in{\bf
N}$ and integer $k$ such that $0\leq k\leq n$: \be \label{2.5}
P_{x,q,n}({\bf 1}(^n_k))=x^k(x;q)_{n-k}\gamma(^n_k)\,, \ee where
$\gamma(^n_k)$ is a polynomial in $q$. We define also
$\gamma(^0_0)=1$. We will prove in the following theorem that
polynomial $\gamma(^n_k)$ coincides with the Gaussian polynomial
$\displaystyle\gc{n}{k}$.
\medskip

{\bf Theorem 2.1.} {\it Let $n, k\in{\bf N}_0$, $0\leq k\leq n$
and $\gamma(^n_k)$ be a polynomial in $q$ defined by
$(\ref{2.5})$. Then
$$
\gamma(^n_k)=\gc{n}{k}.
$$
}

{\bf Proof.} We introduce the following notation: \be \label{2.6}
p(^n_k):=p_{x,q}(^n_k):= P_{x,q,n}({\bf 1}(^n_k))
=\gamma(^n_k)x^k(x;q)_{n-k}\,. \ee

First, we prove that the following recurrence relation is valid
\be \label{2.7}
p(^n_k)=p(^{n-1}_k)(1-q^{n-k-1}x)+p(^{n-1}_{k-1})q^{n-k}x\,. \ee
We apply the formula of total probability. The events ${\bf
1}(^{n-1}_j)$ ($j=0, 1, \ldots , n-1$) are disjoint with union
$\Omega_n$. Since $P({\bf 1}(^n_k)|\,\,{\bf 1}(^{n-1}_j))=0$ for
all $j$ but $j=k-1$ and $j=k$, we may write
\begin{eqnarray*}
p(^n_k)&=&P({\bf 1}(^n_k))= P({\bf 1}(^n_k)|\,{\bf
1}(^{n-1}_k))P({\bf 1}(^{n-1}_k))+
P({\bf 1}(^n_k)|\,{\bf 1}(^{n-1}_{k-1}))P({\bf 1}(^{n-1}_{k-1}))\\
&=&P({\bf 0}^{(n)}|\,{\bf 1}(^{n-1}_k))P({\bf 1}(^{n-1}_k))+
P({\bf 1}^{(n)}|\,{\bf 1}(^{n-1}_{k-1}))P({\bf 1}(^{n-1}_{k-1}))\\
&=&(1-q^{n-k-1}x)p(^{n-1}_k)+q^{n-k}xp(^{n-1}_{k-1})\,.
\end{eqnarray*}

It follows from (\ref{2.7}) that polynomials $\gamma(^n_k)$
satisfy the following recurrence  identity: \be \label{2.8}
\gamma(^n_k)=\gamma(^{n-1}_k)+q^{n-k}\gamma(^{n-1}_{k-1}). \ee
Indeed, if we insert (\ref{2.6}) into (\ref{2.7}) and use
(\ref{1.2}), we obtain
\begin{eqnarray*}
\gamma(^n_k)x^k(x;q)_{n-k}&=&xq^{n-k}\gamma(^{n-1}_{k-1})x^{k-1}
(x;q)_{(n-1)-(k-1)}\\
&&+(1-q^{n-k-1}x)\gamma(^{n-1}_k)x^k(x;q)_{n-1-k}\\
&=&q^{n-k}\gamma(^{n-1}_{k-1})x^{k}
(x;q)_{n-k}+\gamma(^{n-1}_k)x^k(x;q)_{n-k}\,.
\end{eqnarray*}
Dividing by $x^k(x;q)_{n-k}$, we conclude that (\ref{2.8}) is
true.

Next, (\ref{2.5}) yields that $\gamma(^n_0)=1$ and
$\gamma(^n_n)=1$ for all $n\in{\bf N}$. Therefore polynomials
$\gamma(^i_j)$ and Gaussian polynomials satisfy the same
recurrence identity and boundary conditions. This completes the
proof of Theorem~2.1. $\Box$
\medskip

Thus, (\ref{2.5}) can be written as \be \label{2.9}
p_{x,q}(^n_k)=P({\bf 1}(^n_k))=x^k(x;q)_{n-k}\gc{n}{k}\,. \ee The
next theorem contains new proofs of some known facts (see, for
example, [1], formulas (3.3.4), (3.3.9), (3.3.10)). Our proofs are
based upon probabilistic ideas and are very simple.
\medskip

{\bf Theorem 2.2.} 1) {\it If $n, k\in{\bf N}_0$, $0\leq k\leq n$,
then \be \label{2.10} \gc{n}{k}=\gc{n-1}{k}q^k+\gc{n-1}{k-1}\,.
\ee

$2)$ If $n\in{\bf N}$, then
\begin{eqnarray}
\label{2.11} \sum\limits_{k=0}^n(-1)^k\gc{n}{k}q^{{n-k\choose
2}}=0\,.
\end{eqnarray}

$3)$ If $a, b, m\in{\bf N}_0$, then \be \label{2.12}
\gc{a+b}{m}={1\over 2}\sum \limits_{j,k\geq 0\atop
j+k=m}\gc{a}{j}\gc{b}{k} q^{-jk}\left(q^{ak}+q^{bj}\right)\,. \ee

$4)$ If $a, b, c, m\in{\bf N}_0$, then \begin{eqnarray}
\label{2.13} \gc{a+b+c}{m}={1\over 6}\sum \limits_{j,k,l\geq
0\atop j+k+l=m}&&\!\!\!\!\!\!\gc{a}{j}\gc{b}{k}\gc{c}{l} q^{-(jk+kl+lj)}\nonumber\\
&\times&(q^{a(m-j)+bl}+q^{a(m-j)+ck}+q^{b(m-k)+al} \nonumber\\
&+& q^{b(m-k)+cj}+q^{c(m-l)+bj}+q^{c(m-l)+ak})\,.
\end{eqnarray}

$5)$ If $m,n\in{\bf N}_0$, then \be \label{2.14}
\gc{n+m+1}{m+1}=\sum\limits_{j=0}^nq^j\gc{m+j}{m}\,. \ee }
\medskip

{\bf Proof.} 1) First, let us prove that probabilities $p(^n_k)$
satisfy the following recurrence identity \be \label{2.15}
p_{x,q}(^n_k)=p_{qx,q}(^{n-1}_k)(1-x)+p_{x,q}(^{n-1}_{k-1})x\,.
\ee (See definition of $p_{x,q}(^n_k)$ in (\ref{2.9}).) The events
${\bf 0}^{(1)}$ and ${\bf 1}^{(1)}$ form a partition of the sample
space $\Omega_n$. Hence, by the formula of total probability,
 by (\ref{2.1}), and (\ref{2.2}), we obtain
\begin{eqnarray*}
p_{x,q}(^n_k)&=&P_{x,q,n}({\bf 1}(^n_k))\\
&=&P_{x,q,n}({\bf 1}(^n_k)|\,{\bf 0}^{(1)})P_{x,q,n}({\bf
0}^{(1)})+
P_{x,q,n}({\bf 1}(^n_k)|\,{\bf 1}^{(1)})P_{x,q,n}({\bf 1}^{(1)})\\
&=&P_{qx,q,n-1}({\bf 1}(^{n-1}_k))(1-x)+
P_{x,q,n-1}({\bf 1}(^{n-1}_{k-1}))x\\
&=&p_{qx,q}(^{n-1}_k)(1-x)+p_{x,q}(^{n-1}_{k-1})x\,.
\end{eqnarray*}

Inserting (\ref{2.9}) into (\ref{2.15}), we get
\begin{eqnarray*}
\gc{n}{k}x^k(x;q)_{n-k}&=&(1-x)\gc{n-1}{k}q^kx^k(qx;q)_{n-k-1}\\
&&+ x\gc{n-1}{k-1}x^{k-1}(x;q)_{(n-1)-(k-1)}\,.
\end{eqnarray*}
Using (\ref{1.2}) and dividing by $x^k(x;q)_{n-k}$, we obtain
(\ref{2.10}).

2) The events ${\bf 1}(^n_k)$ ($k=0,1,2,\ldots ,n$) form a
partition of the sample space $\Omega_n$. Therefore
$$
\sum\limits_{k=0}^np(^n_k)=1.
$$
Inserting (\ref{2.9}) into this formula and equating to zero
coefficient of $x^n$ in the left-hand side, we obtain
(\ref{2.11}).

Analogously, one can obtain identities for $q$-binomial
coefficients equating to zero coefficients of $x^2$,\ldots,
$x^{n-1}$.

3) The events ${\bf 1}(^a_0), {\bf 1}(^a_1), {\bf 1}(^a_2), \ldots
,{\bf 1}(^a_a)$ form a partition of the sample space $\Omega_n$.
The formula of total probability yields \be \label{2.16} P({\bf
1}(^{a+b}_m))= \sum\limits_{j=0}^aP({\bf 1}(^{a+b}_m)|\,{\bf
1}(^a_j))P({\bf 1}(^a_j)). \ee We calculate the conditional
probabilities at the right-hand side of (\ref{2.16}). Obviously,
$P({\bf 1}(^{a+b}_m)|\,{\bf 1}(^a_j))=0$, if $j>m$. For $j\leq m$
we have by (\ref{2.1}), (\ref{2.50}), and (\ref{2.2}):
$$
P_{x,q,a+b}({\bf 1}(^{a+b}_m)|\,{\bf 1}(^a_j))=
P_{xq^{a-j},q,b}({\bf 1}(^b_{m-j}))\,.
$$
Thus, by (\ref{1.2}), the $j^{th}$ term at the right-hand side of
(\ref{2.16}) equals
\begin{eqnarray}
\label{2.17} &&\gc{b}{m-j}x^{m-j}q^{(a-j)(m-j)}
(xq^{a-j};q)_{b-m+j}\cdot \gc{a}{j}x^j
(x;q)_{a-j}=\nonumber\\
&&=\gc{b}{m-j}\gc{a}{j}x^mq^{(a-j)(m-j)}(x;q)_{a+b-m}\,.
\end{eqnarray}
After substituting (\ref{2.17}) into (\ref{2.16}) and
 division by $x^m(x;q)_{a+b-m}$, we get
\begin{eqnarray*}
\gc{a+b}{m}=
\sum\limits_{j=0}^m\gc{b}{m-j}\gc{a}{j}q^{(a-j)(m-j)}\,.
\end{eqnarray*}
Replacing $m-j$ with $k$ shows that
\begin{eqnarray}
\label{2.19} \gc{a+b}{m}=\sum\limits_{j,k\geq 0\atop j+k=m}
\gc{a}{j}\gc{b}{k}q^{(a-j)k}\,.
\end{eqnarray}
If we first replace in (\ref{2.19}) $a$ with $b$ and $b$ with $a$,
and then $j$ with $k$ and $k$ with $j$, we get
\begin{eqnarray}
\label{2.20} \gc{a+b}{m}=\sum\limits_{j,k\geq 0\atop j+k=m}
\gc{a}{j}\gc{b}{k}q^{j(b-k)}\,.
\end{eqnarray}
Summing (\ref{2.19}) and (\ref{2.20}), we get (\ref{2.12}).

4) Let us define the events $B^{a,b}_{j,k}$ as follows:
\begin{eqnarray*}
 B^{a,b}_{j,k}:=\{(\varepsilon_1,\ldots
,\varepsilon_{a+b})\in\Omega_{a+b}:
&\!\!\!\!\!& \#\{l:1\leq l\leq a\,, \varepsilon_l=1\}=j\,,\\
&\!\!\!\!\!&\#\{l:a+1\leq l\leq a+b\,, \varepsilon_l=1\}=k\}\,.
\end{eqnarray*}
The events $B^{a,b}_{j,k}$ ($0\leq j\leq a$, $0\leq k\leq b$,
$m-c\leq j+k\leq m$) form a partition of $\Omega_{a+b}$.
Therefore, by the formula of total probability \be \label{2.21}
P({\bf 1}(^{a+b+c}_m))= \sum\limits_ {{0\leq j\leq a\atop 0\leq
k\leq b}\atop m-c\leq j+k\leq m}P({\bf
1}(^{a+b+c}_m)|\,B^{a,b}_{j,k}) P (B^{a,b}_{j,k})\,, \ee where
$P=P_{x,q,a+b+c}$. By the  definition of the probability $p$ on
$\Omega_n$ (see (\ref{2.1}), (\ref{2.50}), (\ref{2.2})) we have
\begin{eqnarray*}
&&P_{x,q,a+b+c}({\bf 1}(^{a+b+c}_m)|\,B^{a,b}_{j,k})=
P_{q^{a+b-j-k}x,q,c}({\bf 1}(^{c}_{m-j-k}))=\\
&&=\gc{c}{m-j-k}q^{(a+b-j-k)(m-j-k)}x^{m-j-k}(q^{a+b-j-k}x,q)_{j+k-(m-c)}\,.
\end{eqnarray*}
Using the multiplication theorem of probability and (\ref{1.2}),
we get
\begin{eqnarray*}
&&P_{x,q,a+b+c}(B^{a,b}_{j,k})= P_{x,q,a+b+c}(B^{a,b}_{j,k}|\,{\bf
1}(^a_j))
P_{x,q,a+b+c}({\bf 1}(^a_j))\\
&&=P_{x,q,a+b}({\bf 1}(^b_k)|\,{\bf 1}(^a_j))P_{x,q,a}({\bf
1}(^a_j))= P_{q^{a-j}x,q,b}({\bf 1}(^b_k))P_{x,q,a}({\bf
1}(^a_j))\\
&&=\gc{b}{k}(q^{a-j}x)^k(q^{a-j}x,q)_{b-k}\cdot
\gc{a}{j}x^j(x;q)_{a-j}\\
&&=\gc{a}{j}\gc{b}{k}q^{k(a-j)}x^{j+k}(x;q)_{a+b-j-k}\,.
\end{eqnarray*}
Thus, the $(j,k)^{th}$ term at the right-hand side of (\ref{2.21})
equals
\begin{eqnarray}
\label{2.22}
\gc{a}{j}\gc{b}{k}\gc{c}{m-j-k}x^mq^{(a+b-j-k)(m-j-k)+k(a-j)}
(x;q)_{a+b+c-m}\,.
\end{eqnarray}
After substituting (\ref{2.22}) into (\ref{2.21}), replacing
$m-j-k$ by $l$, and dividing the result by  $x^m(x;q)_{a+b+c-m}$,
we obtain
\begin{eqnarray}
\label{2.71} \gc{a+b+c}{m}=\sum \limits_{j,k,l\geq 0\atop
j+k+l=m}\gc{a}{j}\gc{b}{k}\gc{c}{l} q^{-(jk+kl+lj)}
q^{a(m-j)+bl}\,.
\end{eqnarray}
Formula (\ref{2.13}) follows from (\ref{2.71}) with the help of
suitable permutations of $a, b, c$ and $j, k, l$ and summing. (See
the end of the proof of (\ref{2.12}).)

5) Let us calculate the probability $P({\bf 1}(^{n+m+1}_{m+1}))$
in two ways. We have by (\ref{2.9}) \be \label{2.23}
P_{x,q,n+m+1}({\bf 1}(^{n+m+1}_{m+1}))=
\gc{n+m+1}{m+1}x^{m+1}(x;q)_n\,. \ee Now we calculate this
probability in a different way. For $j=0, 1, 2, \ldots , n$, we
introduce the following events
\begin{eqnarray*}
C_j:={\bf 1}(^{n+m+1}_{m+1})\cap
\{&\!\!\!\!\!\!&(\varepsilon_1,\ldots,\varepsilon_{n+m+1})\in\Omega_{n+m+1}:\\
&\!\!\!\!\!\!&\varepsilon_{m+j+1}=1, \varepsilon_{m+j+2}=\ldots
=\varepsilon_{n+m+1}=0\} \,.
\end{eqnarray*}
(We distinguished the last outcome $1$ in the sequence of $n+m+1$
trials ${\cal E}$.) Any element $\omega$ of $C_j$ can be written
as $\omega=(\omega^{\prime},1,0_{(n-j)})$, where
$\omega^{\prime}\in{\bf 1}(^{m+j}_m)$. The events $C_j$ are
mutually exclusive and $\cup_{j=0}^nC_j={\bf 1}(^{n+m+1}_{m+1})$.
Therefore: \be \label{2.24} P ({\bf 1}(^{n+m+1}_{m+1}))=
\sum\limits_{j=0}^nP (C_j)\,. \ee Using (\ref{2.9}) and
(\ref{1.2}), we obtain:
\begin{eqnarray}
\label{2.25}
&&P_{x,q,n+m+1}(C_j)=\nonumber\\
&&=P_{x,q,m+j}({\bf
1}(^{m+j}_{m}))\cdot(xq^j)\cdot(1-xq^j)(1-xq^{j+1})\cdots
(1-xq^{n-1})\nonumber\\
&&=\gc{m+j}{m}x^m(x;q)_j\cdot xq^j\cdot(xq^j;q)_{n-j}\nonumber\\
&&=q^j\gc{m+j}{m}x^{m+1}(x;q)_n\,.
\end{eqnarray}
Substituting (\ref{2.23}) and (\ref{2.25}) into (\ref{2.24}) and
dividing by $x^{m+1}(x;q)_n$, we obtain (\ref{2.14}). $\Box$
\medskip

We may distinguish the first rather than the last - or the first
and the last - outcome $1$ in the sequence of trials. This leads
to the following theorem.
\medskip

{\bf Theorem 2.3.} {\it If $m,n\in{\bf N}_0$, then \be
\label{2.26}
\gc{n+m+1}{m+1}=\sum\limits_{j=0}^nq^{(m+1)j}\gc{n+m-j}{m}\,, \ee
\be \label{2.27} \gc{n+m+2}{m+2}={1\over 2}\sum\limits_{j,k\geq
0\atop j+k\leq n}
\gc{n+m-j-k}{m}q^{n-(j+k)}\left(q^{(m+2)j}+q^{(m+2)k}\right)\,.
\ee }

{\bf Proof.} 1) We prove (\ref{2.26}). We calculate $P({\bf
1}(^{n+m+1}_{m+1}))$ in two ways. Evidently, we have (\ref{2.23}).
Alternatively, for  $j=0, 1, 2, \ldots, n$, we introduce the
following events
\begin{eqnarray*}
D_j:={\bf 1}(^{n+m+1}_{m+1})\cap
\{(\varepsilon_1,\ldots,\varepsilon_{n+m+1})\in\Omega_{n+m+1}:
\varepsilon_1=\ldots =\varepsilon_j=0, \varepsilon_{j+1}=1\} \,.
\end{eqnarray*}
$D_j$ is formed by elements $\omega=(0_{(j)},1,\omega^{\prime})$
where $\omega^{\prime}\in{\bf 1}(^{n+m+1}_m)$.  Therefore
\begin{eqnarray}
\label{2.28} P(D_j)&=&(x;q)_j \cdot
(xq^{j})\cdot P_{xq^j,q,n+m+-j}({\bf 1}(^{n+m-j}_m))\nonumber\\
&=&(x;q)_j\cdot\gc{n-j+m}{m}\cdot(xq^{j})^m(xq^{j};q)_{n-j}
\nonumber\\
&=&\gc{n-j+m}{m}x^{m+1}(x;q)_nq^{(m+1)j}\,.
\end{eqnarray}
The events $D_j$ are mutually exclusive and
$\cup_{j=0}^{n}D_j={\bf 1}(^{n+m+1}_{m+1})$. Therefore, by the
additive property of probability, $ P({\bf
1}(^{n+m+1}_{m+1}))=\sum_{j=0}^{n}P(D_j)$. Inserting (\ref{2.23})
and (\ref{2.28}) into this formula and dividing by
$x^{m+1}(x;q)_n$, we get (\ref{2.26}).

2)  We prove (\ref{2.27}). For $j\geq 0$ i $k\geq 0$ such that
$j+k\leq n$, we consider events
\begin{eqnarray*}
E_{j,k}:={\bf 1}(^{n+m+2}_{m+2})\cap
\{&\!\!\!\!\!\!&(\varepsilon_1,\ldots,\varepsilon_{n+m+2})\in\Omega_{n+m+2}:\\
&\!\!\!\!\!\!&\varepsilon_1=\ldots\varepsilon_j=0, \varepsilon_{j+1}=1,\\
&\!\!\!\!\!\!&\varepsilon_{n+m+3-k}=1,
\varepsilon_{n+m+4-k}=\ldots = \varepsilon_{n+m+2}=0\}\,.
\end{eqnarray*}
The event $E_{j,k}$ is formed by elements
$\omega=(0_{(j)},1,\omega^{\prime},1,0_{(k)})$ where
$\omega^{\prime}\in{\bf 1}(^{n+m-j-k}_m)$. We calculate the
probability of the events $E_{j,k}$ using (\ref{2.1}),
(\ref{2.50}), (\ref{2.2}):
\begin{eqnarray*}
&&P_{x,q,n+m+2}(E_{j,k})=\\
&&=(x;q)_j \cdot xq^{j} \cdot P_{q^jx,q,n+m-j-k}({\bf
1}(^{n+m-j-k}_{m})) \cdot xq^{j+(n-j-k)}\cdot (xq^{n-k};q)_k \,.
\end{eqnarray*}
Since
$$
P_{q^jx,q,n+m-j-k}({\bf 1}(^{n+m-j-k}_{m}))=
\gc{n+m-j-k}{m}(xq^{j})^{m} (xq^j,q)_{n-(j+k)}\,,
$$
it follows by (\ref{1.2}) that
\begin{eqnarray}
\label{2.29} P(E_{j,k})=
\gc{n+m-j-k}{m}q^{mj+j+n-k}x^{m+2}(x;q)_n\,.
\end{eqnarray}

Since $E_{j,k}$ are mutually exclusive and ${\bf
1}(^{n+m+2}_{m+2}) =\cup_{j,k}E_{j,k}$, the following equality
holds: $P({\bf 1}(^{n+m+2}_{m+2}))= \sum\limits_{j,k}P(E_{j,k})$.
Substituting (\ref{2.29}) into this equality, using (\ref{2.9}),
and dividing by $x^{m+2}(x;q)_n$, we get \be \label{2.30}
\gc{n+m+2}{m+2}=\sum\limits_{j,k\geq 0\atop j+k\leq n}
\gc{n+m-j-k}{m}q^{mj+n+j-k}\,. \ee Replacing  $j$ in (\ref{2.30})
by $k$, and $k$ by $j$, we obtain \be \label{2.31}
\gc{n+m+2}{m+2}=\sum\limits_{j,k\geq 0\atop j+k\leq n}
\gc{n+m-j-k}{m}q^{mk+n+k-j}\,. \ee Summing (\ref{2.30}) and
(\ref{2.31}) we get (\ref{2.27}). $\Box$
\bigskip

\section{$q$-polynomial coefficients}
\setcounter{equation}{0}
\bigskip

Let $n\in{\bf N}$, $u_1>0$, $u_2>0$, $U:=u_1+u_2<1$, $0<q<1$. Let
${\cal E}$ be a trial with four outcomes $0,1,2,\ast$. We consider
the following random process. The trial ${\cal E}$ is repeated $n$
times. We assume that the probabilities of the outcomes $0,1,2$ in
the first trial are equal to \be \label{3.1} 1-u_1-u_2=1-U\,,\quad
u_1\,,\quad u_2\,, \ee respectively. Consequently, the probability
of the outcome $\ast$ in the first trial is equal to zero.

Let $m$ be a positive integer, $m<n$. Suppose that in the first
$m$ trials the outcomes $0, 1, 2$ have happened $i_0, i_1, i_2$
times, respectively, with $i_0+i_1+i_2=m$ (so, there is no outcome
$\ast$ in the first $m$ trials).  Then we assume the probabilities
of outcomes $0, 1, 2$ in the $(m+1)^{th}$ trial to be equal \be
\label{3.2} 1-q^{i_0}U\,,\quad q^{i_0}u_1\,,\quad
q^{i_0+i_1}u_2\,, \ee respectively. Consequently, the outcome
$\ast$ happens in the $(m+1)^{th}$ trial with the probability \be
\label{3.3} 1-q^{i_0}u_1-q^{i_0+i_1}u_2-
(1-q^{i_0}u_1-q^{i_0}u_2)= q^{i_0}u_2(1-q^{i_1})\,. \ee We assume
also that if the outcome $\ast$ happens in the $k^{th}$ trial,
then $\ast$ will happen in the $(k+1)^{th}$ trial with the
probability $1$.

We construct a probability space corresponding to the random
process described above. The sample space $\Omega_n$ consists of
all sequences
$\omega=(\varepsilon_1,\varepsilon_2,\ldots,\varepsilon_n)$ of the
length $n$, such that its elements $\varepsilon_j$ are equal to
$0,1,2,\ast$, and the following condition is valid: if
$\varepsilon_k=\ast$ for some $k$ ($1<k<n$), then
$\varepsilon_l=\ast$ for any $l=k+1, k+2, \ldots , n$. We define
the probability $p_{u_1,u_2,U,q,n}(\omega)$ of the elementary
event $\omega=(\varepsilon_1,\varepsilon_2,\ldots,\varepsilon_n)$
as \be \label{3.4} p_{u_1,u_2,U,q,n}
((\varepsilon_1,\varepsilon_2,\ldots,\varepsilon_n))= f_1\cdot
f_2\cdot\ldots\cdot f_n\,, \ee where, according to (\ref{3.1}),
\begin{eqnarray}
\label{3.5} f_1=\cases{1-U, & if $\varepsilon_1=0$,\cr u_1, & if
$\varepsilon_1=1$,\cr u_2, & if $\varepsilon_1=2$,\cr 0, & if
$\varepsilon_1=\ast$,\cr}
\end{eqnarray}
and according to (\ref{3.2}) and (\ref{3.3}), for any $m$, $1\leq
m\leq n-1$,
\begin{eqnarray}
\label{3.6} f_{m+1}=\cases{1-q^{i_0}U\,,& if
$\varepsilon_{m+1}=0$\,,\cr q^{i_0}u_1\,, & if
$\varepsilon_{m+1}=1$\,,\cr q^{i_0+i_1}u_2\,,& if
$\varepsilon_{m+1}=2$\,,\cr q^{i_0}u_2(1-q^{i_1})\,, & if
$\varepsilon_{m+1}=\ast$\,,\cr}
\end{eqnarray}
if $i_0+i_1+i_2=m$ and
\begin{eqnarray*}
&&\#\{j:1\leq j\leq m, \varepsilon_j=0\}=i_0\,,\\
&&\#\{k:1\leq k\leq m, \varepsilon_k=1\}=i_1\,,\\
&&\#\{l:1\leq l\leq m, \varepsilon_l=2\}=i_2\,.
\end{eqnarray*}
We define also \be \label{3.7} \varepsilon_{m+1}=\ast\quad\mbox{
and }\quad f_{m+1}=1,\qquad\mbox{ if }\varepsilon_m=\ast \ee (see
the end of the second paragraph of this section). We define the
probability of an event $A\subset\Omega_n$ as follows \be
\label{3.8} P_{u_1,u_2,U,q,n}(A):= \sum\limits_{\omega\in
A}p_{u_1,u_2,U,q,n}(\omega)\,. \ee Although $U=u_1+u_2$, it is
useful to write $U$ as a subscript of the probability sign $P$.
Sometimes we write $P$ instead of $ P_{u_1,u_2,U,q,n}$. As in
Section~2, we see that
$$
P_{u_1,u_2,U,q,m}(A)=P_{u_1,u_2,U,q,n}(A)\,,
$$
if $m<n$ and the event $A\subset\Omega_n$ does not depend on
repetitions of the trial ${\cal E}$ with the numbers $m+1,
m+2,\ldots , n$.

Let us introduce notation for some events of the sample space
$\Omega_n$. For $i\in\{0,1,2,\ast\}$ and integer $m$, $1\leq m\leq
n$, and $k$, $0\leq k\leq n$, we use definitions (\ref{2.3}),
(\ref{2.4}) from Section~2. For $m\in{\bf N}$, $j\in\{0,1,2\}$,
$i_0,i_1,i_2\in{\bf N}_0$ such that $i_0+i_1+i_2=m$, we define
\begin{eqnarray*}
\af{A}{m}{i_0}{i_1}{i_2} := {\bf 0}(^m_{i_0})\cap{\bf
1}(^m_{i_1})\cap{\bf 2}(^m_{i_2})\,.
\end{eqnarray*}
The event $\af{A}{m}{i_0}{i_1}{i_2}$ can be described as follows:
outcomes $0, 1, 2$ happen $i_0, i_1, i_2$ times, respectively, in
the repetitions of the trial ${\cal E}$ with the numbers $1,
2,\ldots , m$. We set $\af{A}{n}{i_0}{i_1}{i_2}=\emptyset$ if
$i_j<0$ or $i_j>m$ for some $j\in\{0, 1, 2 \}$.

We may write (\ref{3.1}), (\ref{3.2}), (\ref{3.3}) in the
introduced notation as follows ($P=P_{u_1,u_2,U,q,n}$): \be
\begin{array}{rclcrcl}
\label{3.9} P({\bf 0}^{(1)})&=&1-u_1-u_2\,,&\quad&
P({\bf 1}^{(1)})&=&u_1\,,\\
P({\bf 2}^{(1)})&=&u_2\,,&\quad& P({\bf\ast}^{(1)})&=&0\,;
\end{array}
\ee if $m\in{\bf N}$, $m<n$, $i_0+i_1+i_2=m$, then \be
\label{3.10}
\begin{array}{rcl}
P\left({\bf 0}^{(m+1)}|
\af{A}{m}{i_0}{i_1}{i_2}\right)&=&1-q^{i_0}U\,,\\
P\left({\bf 1}^{(m+1)}|
\af{A}{m}{i_0}{i_1}{i_2}\right)&=&q^{i_0}u_1\,,\\
P\left({\bf 2}^{(m+1)}|
\af{A}{m}{i_0}{i_1}{i_2}\right)&=&q^{i_0+i_1}u_2\,,\\
P\left({\bf\ast}^{(m+1)}|\af{A}{m}{i_0}{i_1}{i_2}\right)&=&
q^{i_0}(1-q^{i_1})u_2\,.
\end{array}
\ee Formula (\ref{3.7}) means that
\begin{eqnarray*}
P({\bf j}^{(l+1)}|\ast^{(l)})=0\quad\mbox{ for } j=0,1,2\,;\qquad
P({\bf\ast}^{(l+1)}|\ast^{(l)})=1
\end{eqnarray*}
for all $l=1, 2, \ldots, n-1$.

We give now the probabilities of some elementary events
$\omega\in\af{A}{n}{i_0}{i_1}{i_2}$:
\begin{eqnarray*}
&&p((1_{(i_1)},2_{(i_2)}, 0_{(i_0)}))=
u_1^{i_1}u_2^{i_2}(U;q)_{i_0}q^{i_1i_2}\,,\\
&&p((2_{(i_2)},1_{(i_1)}, 0_{(i_0)}))=
u_1^{i_1}u_2^{i_2}(U;q)_{i_0}\,,\\
&&p((0_{(i_0)},1_{(i_1)}, 2_{(i_2)}))=
u_1^{i_1}u_2^{i_2}(U;q)_{i_0}q^{i_0u_1+(i_0+i_1)i_2}\,,\\
&&p((0,1,2,0,1,2_{(i_2-1)}, 0_{(i_0-2)}, 1_{(i_1-2)}))=
u_1^{i_1}u_2^{i_2}(U;q)_{i_0}q^{i_0(i_1-2)+4i_2+1}\,.
\end{eqnarray*}
It can be readily  seen that the probability of every elementary
event $\omega\in\af{A}{n}{i_0}{i_1}{i_2}$ $(i_0+i_1+i_2=n)$ equals
$u_1^{i_1}u_2^{i_2}(U;q)_{i_0}q^r$, where $r$ is an integer such
that $0\leq r\leq i_0i_1+i_0i_2+i_1i_2$. Therefore, \be
\label{3.11}
\af{p}{n}{i_0}{i_1}{i_2}:=P\left(\af{A}{n}{i_0}{i_1}{i_2}\right)=
u_1^{i_1}u_2^{i_2}(U;q)_{i_0}\af{c}{n}{i_0}{i_1}{i_2}\,, \ee where
$\af{c}{n}{i_0}{i_1}{i_2}=c(n,i_0,i_1,i_2;q)$ is a polynomial in
$q$. We set $\af{c}{0}{0}{0}{0}=1$ and
$\af{c}{n}{i_0}{i_1}{i_2}=0$ if one of the numbers $i_0, i_1, i_2$
is negative.

The following theorem is an analogue of Theorem~2.1. It shows that
the polynomial $\af{c}{n}{i_0}{i_1}{i_2}$ coincides with the
$q$-polynomial coefficient $\displaystyle\qp{n}{i_0}{i_1}{i_2}$.
It is a key theorem of this section.
\medskip

{\bf Theorem 3.1.} {\it Let $n\in{\bf N}_0$, $i_0, i_1, i_2\in{\bf
N}_0$, $i_0+i_1+i_2=n$, and $\af{c}{n}{i_0}{i_1}{i_2}$ be a
polynomial in $q$ defined by $(\ref{3.11})$. Then
$$
\af{c}{n}{i_0}{i_1}{i_2}=\left[n \atop i_0,i_1,i_2\right]_q\,.
$$
}

{\bf Proof.} First we prove a recurrence relation for coefficients
$\af{c}{n}{i_0}{i_1}{i_2}$. We show that if $n\geq 1$, $i_0, i_1,
i_2\geq 0$, $i_0+i_1+i_2=n$, then
\begin{eqnarray}
\label{3.12} \af{c}{n}{i_0}{i_1}{i_2}=
\af{c}{n-1}{i_0-1}{i_1}{i_2}+ q^{i_0}\af{c}{n-1}{i_0}{i_1-1}{i_2}+
q^{i_0+i_1}\af{c}{n-1}{i_0}{i_1}{i_2-1}\,.
\end{eqnarray}

Since the events $\af{A}{n-1}{j_0}{j_1}{j_2}$ ($j_0, j_1, j_2\geq
0$, $j_0+j_1+j_2=n-1$) and
$B_{n-1}:=\{(\varepsilon_1,\ldots,\varepsilon_n)
\in\Omega_n:\varepsilon_k=\ast\,\,\,(\exists k, 1\leq k\leq
n-1)\}$ form a partition of the sample space $\Omega_n$ and since
$P\left(\af{A}{n}{i_0}{i_1}{i_2}|B_{n-1}\right)=0$, it follows by
the formula of total probability that
\begin{eqnarray}
\label{3.13} &&\af{p}{n}{i_0}{i_1}{i_2}=
P\left(\af{A}{n}{i_0}{i_1}{i_2}\right)=\nonumber\\
&&=P\left(\af{A}{n}{i_0}{i_1}{i_2}|\af{A}{n-1}{i_0-1}{i_1}{i_2}\right)
P\left(\af{A}{n-1}{i_0-1}{i_1}{i_2}\right)\nonumber\\
&&+P\left(\af{A}{n}{i_0}{i_1}{i_2}|\af{A}{n-1}{i_0}{i_1-1}{i_2}\right)
P\left(\af{A}{n-1}{i_0}{i_1-1}{i_2}\right)\nonumber\\
&&+P\left(\af{A}{n}{i_0}{i_1}{i_2}|\af{A}{n-1}{i_0}{i_1}{i_2-1}\right)
P\left(\af{A}{n-1}{i_0}{i_1}{i_2-1}\right)\nonumber\\
&&=P\left({\bf 0}^{(n)}|\af{A}{n-1}{i_0-1}{i_1}{i_2}\right)
P\left(\af{A}{n-1}{i_0-1}{i_1}{i_2}\right)\nonumber\\
&&+P\left({\bf 1}^{(n)}|\af{A}{n-1}{i_0}{i_1-1}{i_2}\right)
P\left(\af{A}{n-1}{i_0}{i_1-1}{i_2}\right)\nonumber\\
&&+P\left({\bf 2}^{(n)}|\af{A}{n-1}{i_0}{i_1}{i_2-1}\right)
P\left(\af{A}{n-1}{i_0}{i_1}{i_2-1}\right)\nonumber\\
&&=(1-q^{i_0-1}U)\af{p}{n-1}{i_0-1}{i_1}{i_2}+
q^{i_0}u_1\af{p}{n-1}{i_0}{i_1-1}{i_2}+
q^{i_0+i_1}u_2\af{p}{n-1}{i_0}{i_1}{i_2-1}\,.
\end{eqnarray}
Substituting (\ref{3.11}) into (\ref{3.13}), using the formula
$(1-q^{i_0-1}U)(U;q)_{i_0-1}=(U;q)_{i_0}$, and dividing by
$u_1u_2(U;q)_{i_0}$, we obtain (\ref{3.12}).

Now we prove the result by induction on $n$. By definition
$\displaystyle\af{c}{0}{0}{0}{0}=\qp{0}{0}{0}{0}$. We assume that
the equality
$\displaystyle\af{c}{n-1}{j_0}{j_1}{j_2}=\qp{n-1}{j_0}{j_1}{j_2}$
holds for all $j_0, j_1, j_2\geq 0$ such that $j_0+j_1+j_2=n-1$.
We will show that
$\displaystyle\af{c}{n}{i_0}{i_1}{i_2}=\qp{n}{i_0}{i_1}{i_2}$ for
all $i_0, i_1, i_2\geq 0$ such that $i_0+i_1+i_2=n$. Applying
(\ref{3.12}) we get:
\begin{eqnarray*}
&&\af{c}{n}{i_0}{i_1}{i_2}=\af{c}{n-1}{i_0-1}{i_1}{i_2}+
q^{i_0}\af{c}{n-1}{i_0}{i_1-1}{i_2}+
q^{i_0+i_1}\af{c}{n-1}{i_0}{i_1}{i_2-1}\\
&&={(q)_{n-1}\over (q)_{i_0}(q)_{i_1}(q)_{i_2}}
\left((1-q^{i_0})+q^{i_0}(1-q^{i_1})+q^{i_0+i_1}(1-q^{i_2})\right)\\
&&={(q)_{n-1}\over (q)_{i_0}(q)_{i_1}(q)_{i_2}}(1-q^{n})=
{(q)_{n}\over (q)_{i_0}(q)_{i_1}(q)_{i_2}}\,,
\end{eqnarray*}
and the proof of the theorem is complete. $\Box$
\medskip

Therefore, we may write \be \label{3.14}
\af{p}{n}{i_0}{i_1}{i_2}=P\left(\af{A}{n}{i_0}{i_1}{i_2}\right)=
u_1^{i_1}u_2^{i_2}(U;q)_{i_0}\qp{n}{i_0}{i_1}{i_2}\,. \ee

The next theorem contains another recurrence relation for
$q$-polynomial coefficients and an analogue  of the Vandermond
formula.
\medskip

{\bf Theorem 3.2.} 1) {\it If $n, i_0, i_1, i_2\in{\bf N}_0$,
$i_0+i_1+i_2=n$, then
\begin{eqnarray}
\label{3.15} \qp{n}{i_0}{i_1}{i_2}=
q^{i_1+i_2}\qp{n-1}{i_0-1}{i_1}{i_2}+
q^{i_2}\qp{n-1}{i_0}{i_1-1}{i_2}+\qp{n-1}{i_0}{i_1}{i_2-1}\,.
\end{eqnarray}
$2)$  If $n, k\in{\bf N}_0$, $0\leq k\leq n$, $i_0, i_1, i_2\geq
0$, $i_0+i_1+i_2=n$, then
\begin{eqnarray}
\label{3.16}
&&\qp{n}{i_0}{i_1}{i_2}=\\
&&=\sum\limits_{{j_0,j_1,j_2\geq 0\atop
j_0+j_1+j_2=k}}\qp{k}{j_0}{j_1}{j_2}
\qp{n-k}{i_0-j_0}{i_1-j_1}{i_2-j_2}q^{j_0(i_1-j_1)+(j_0+j_1)(i_2-j_2)}\,.\nonumber
\end{eqnarray}
}

{\bf Proof.} 1) Since the events
$\displaystyle\af{A}{1}{1}{0}{0}$,
$\displaystyle\af{A}{1}{0}{1}{0}$,
$\displaystyle\af{A}{1}{0}{0}{1}$ and $B_1:=\{(\ast_{(n)})\}$ form
a partition of the sample space $\Omega_n$ and since $P(B_1)=0$,
we get by the formula of total probability and using (\ref{3.9})
(\ref{3.10}) that
\begin{eqnarray}
\label{3.17} &&\af{p}{n}{i_0}{i_1}{i_2}=\nonumber\\
&&=P\left(\af{A}{n}{i_0}{i_1}{i_2}\right)=
P\left(\af{A}{n}{i_0}{i_1}{i_2}|\af{A}{1}{1}{0}{0}\right)
P\left(\af{A}{1}{1}{0}{0}\right)\nonumber\\
&&+P\left(\af{A}{n}{i_0}{i_1}{i_2}|\af{A}{1}{0}{1}{0}\right)
P\left(\af{A}{1}{0}{1}{0}\right)+
P\left(\af{A}{n}{i_0}{i_1}{i_2}|\af{A}{1}{0}{0}{1}\right)
P\left(\af{A}{1}{0}{0}{1}\right)\nonumber\\
&&=P_{qu_1,qu_2,qU,q,n-1}(\af{A}{n-1}{i_0-1}{i_1}{i_2})(1-U)+
P_{u_1,qu_2,U,q,n-1}(\af{A}{n-1}{i_0}{i_1-1}{i_2})u_1\nonumber\\
&&+\af{p}{n-1}{i_0}{i_1}{i_2-1}u_2\,.
\end{eqnarray}
By the  definition of the probability on $\Omega_n$ (see
(\ref{3.11})), it follows from (\ref{3.17}) that
\begin{eqnarray*}
&&(U;q)_{i_0}u_1^{i_1}u_2^{i_2}\qp{n}{i_0}{i_1}{i_2}=
(1-U)(qU,q)_{i_0-1}(qu_1)^{i_1}(qu_2)^{i_2}
\qp{n-1}{i_0-1}{i_1}{i_2}\\
&&+u_1(U;q)_{i_0}u_1^{i_1-1}(qu_2)^{i_2}\qp{n-1}{i_0}{i_1-1}{i_2}+
u_2(U;q)_{i_0}u_1^{i_1}u_2^{i_2-1}\qp{n-1}{i_0}{i_1}{i_2-1}\,.
\end{eqnarray*}
Using the formula $(1-U)(qU,q)_{i_0-1}=(U;q)_{i_0}$ and dividing
by $u_1^{i_1}u_2^{i_2}(U;q)_{i_0}$, we get (\ref{3.15}).

2) The events $\af{A}{k}{j_0}{j_1}{j_2}$ ($k$ is fixed, $j_0, j_1,
j_2\geq 0$, $j_0+j_1+j_2=k$) and
$B_k:=\{(\varepsilon_1,\ldots,\varepsilon_n)\in\Omega_n:\varepsilon_l=
\ast\,\,\,(\exists l, 1\leq l\leq k)\}$ form a partition of
$\Omega_n$. By the formula of total probability, it follows that
\begin{eqnarray}
\label{3.18}
\af{p}{n}{i_0}{i_1}{i_2}&=&P\left(\af{A}{n}{i_0}{i_1}{i_2}\right)
=P\left(\af{A}{n}{i_0}{i_1}{i_2}|B_k\right)P(B_k)\nonumber\\
&+&\sum\limits_{j_0,j_1,j_2\geq 0\atop j_0+j_1+j_2=k}
P\left(\af{A}{n}{i_0}{i_1}{i_2}|\af{A}{k}{j_0}{j_1}{j_2}\right)
P\left(\af{A}{k}{j_0}{j_1}{j_2}\right)\,.
\end{eqnarray}
Let us first find the conditional probabilities at the right-hand
side of (\ref{3.18}). If the event $\af{A}{k}{j_0}{j_1}{j_2}$ has
occurred, then the probabilities of outcomes $0, 1, 2$ in the
$(k+1)^{th}$ trial ${\cal E}$ are equal to $1-q^{j_0}U$,
$q^{j_0}u_1$, and $q^{j_0+j_1}u_2$, respectively. In $n-k$ trials
with numbers $k+1,k+2,\dots , n$ there must be $i_0-j_0$ outcomes
$0$, $i_1-j_1$ outcomes $1$, and $i_2-j_2$ outcomes $2$.
Therefore:
\begin{eqnarray}
\label{3.19}
&&P\left(\af{A}{n}{i_0}{i_1}{i_2}|\af{A}{k}{j_0}{j_1}{j_2}\right)=\nonumber\\
&&=P_{q^{j_0}u_1,q^{j_0+j_1}u_2,q^{j_0}U,q,n-k}
\left(\af{A}{n-k}{i_0-j_0}{i_1-j_1}{i_2-j_2}\right)\\
&&=\left(q^{j_0}u_1\right)^{i_1-j_1}\left(q^{j_0+j_1}u_2\right)^{i_2-j_2}
(q^{j_0}U,q)_{i_0-j_0}
\qp{n-k}{i_0-j_0}{i_1-j_1}{i_2-j_2}\,.\nonumber
\end{eqnarray}
Substituting (\ref{3.14}) into (\ref{3.18}), using the equality
$P\left(\af{A}{n}{i_0}{i_1}{i_2}|B_k\right)=0$ and (\ref{3.19}),
we get
\begin{eqnarray*}
&&(U;q)_{i_0}u_1^{i_1}u_2^{i_2}\qp{n}{i_0}{i_1}{i_2}=
\sum\limits_{(j_0,j_1,j_2)} \left(q^{j_0}u_1\right)^{i_1-j_1}
\left(q^{j_0+j_1}u_2\right)^{i_2-j_2}\times\\
&&\times(q^{j_0}U,q))_{i_0-j_0}
\qp{n-k}{i_0-j_0}{i_1-j_1}{i_2-j_2} (U;q)_{j_0}u_1^{j_1}u_2^{j_2}
\qp{k}{j_0}{j_1}{j_2}\,.
\end{eqnarray*}
Using (\ref{1.2}) and dividing by $(U;q)_{i_0}u_1^{i_1}u_2^{i_2}$,
we obtain (\ref{3.16}). $\Box$
\medskip

The following theorem contains results that give analogues of
Theorems~2.2 and 2.3.
\medskip

{\bf Theorem 3.3.} {\it For all $i_0$, $i_1, i_2\in{\bf N}_0$ the
following identities are valid:
\begin{eqnarray}
\label{3.20} &&\qp{i_0+i_1+i_2}{i_0}{i_1}{i_2}=
\sum\limits_{k=0}^{i_0}
\qp{k+i_1+i_2}{k}{i_1}{i_2}\left(1-{1-q^k\over
1-q^{k+i_1+i_2}}\right)\,,\\
\label{3.21} &&\qp{i_0+i_1+i_2}{i_0}{i_1}{i_2}=
\left(1-q^{i_1+i_2}\right) \sum\limits_{k=0}^{i_0}
\qp{k+i_1+i_2}{k}{i_1}{i_2}
{q^{(i_1+i_2)(i_0-k)}\over 1-q^{k+i_1+i_2}}\,,\\
\label{3.22}
&&\qp{i_0+i_1+i_2}{i_0}{i_1}{i_2}=\nonumber\\
&&=\left(1-q^{i_1+i_2-1}\right) \sum\limits_{k=0}^{i_0}
\qp{k+i_1+i_2}{k}{i_1}{i_2}
{q^k\left(1-q^{(i_1+i_2)(i_0-k+1)}\right) \over
\left(1-q^{k+i_1+i_2-1}\right) \left(1-q^{k+i_1+i_2}\right)}\,.
\end{eqnarray}
}

{\bf Remark 3.1.} Identities (\ref{2.14}), (\ref{2.26}),
(\ref{2.27}) follow from (\ref{3.20}), (\ref{3.21}), (\ref{3.22}),
respectively, if we take $i_2=0$.

{\bf Proof.} 1) For $r=1,2$ and $j=i_1+i_2, i_1+i_2+1,\ldots,
i_0+i_1+i_2$ we consider the events
$$
C_j^r:=\af{A}{i_0+i_1+i_2}{i_0}{i_1}{i_2}\cap
\{(\varepsilon_1,\ldots,\varepsilon_{i_0+i_1+i_2}):
\varepsilon_j=r,\varepsilon_{j+1}=\ldots =
\varepsilon_{i_0+i_1+i_2}=0\}\,.
$$
By the additive property of probability, we have
\begin{eqnarray}
\label{3.23} P\left(\af{A}{i_0+i_1+i_2}{i_0}{i_1}{i_2}\right)=
\sum\limits_{j=i_1+i_2}^{i_0+i_1+i_2}P\left(C_j^1\right)+
\sum\limits_{j=i_1+i_2}^{i_0+i_1+i_2}P\left(C_j^2\right)\,.
\end{eqnarray}
It is readily seen that \be \label{3.24}
\begin{array}{rl}
&P(C_j^1)= (U;q)_{i_0}u_1^{i_1}u_2^{i_2}q^{j-i_1-i_2}
{\displaystyle\qp{j-1}{j-i_1-i_2}{i_1-1}{i_2}}\,,\\
&P(C_j^2)= (U;q)_{i_0}u_1^{i_1}u_2^{i_2}q^{j-i_2}
{\displaystyle\qp{j-1}{j-i_1-i_2}{i_1}{i_2-1}}\,.
\end{array}
\ee Substituting (\ref{3.14}) and (\ref{3.24}) into (\ref{3.23})
and dividing by $(U;q)_{i_0}u_1^{i_1}u_2^{i_2}$, we get
\begin{eqnarray}
\label{3.25} \qp{i_0+i_1+i_2}{i_0}{i_1}{i_2}&=&
\sum\limits_{j=i_1+i_2}^{i_0+i_1+i_2}
\qp{j-1}{j-i_1-i_2}{i_1-1}{i_2}q^{j-i_1-i_2}+\nonumber\\
&+&\sum\limits_{j=i_1+i_2}^{i_0+i_1+i_2}
\qp{j-1}{j-i_1-i_2}{i_1}{i_2-1}q^{j-i_2}\,.
\end{eqnarray}
Replacing $j$ in (\ref{3.25}) by $k=j-i_1-i_2$ and using the
following identities \be \label{3.26}
\begin{array}{rl}
&{\displaystyle\qp{a+b+c-1}{a}{b-1}{c}=\qp{a+b+c}{a}{b}{c}
{1-q^b\over 1-q^{a+b+c}}}\,,\\
&{\displaystyle\qp{a+b+c-1}{a}{b}{c-1}=\qp{a+b+c}{a}{b}{c}
{1-q^c\over 1-q^{a+b+c}}}\,,
\end{array}
\ee we obtain (\ref{3.20}) after a simple calculation.

2) For $r=1,2$ and $j=0,1,\ldots, i_0$, we consider the events
$$
D_j^r:=\af{A}{i_0+i_1+i_2}{i_0}{i_1}{i_2}\cap
\{(\varepsilon_1,\ldots,\varepsilon_{i_0+i_1+i_2}):
\varepsilon_1=\ldots=\varepsilon_{j}=0, \varepsilon_{j+1}=r\}\,.
$$
Using the additive property of probability, we write
\begin{eqnarray}
\label{3.27} P\left(\af{A}{i_0+i_1+i_2}{i_0}{i_1}{i_2}\right)=
\sum\limits_{j=0}^{i_0}\left(P\left(D_j^1\right)+
P\left(D_j^2\right)\right)\,.
\end{eqnarray}
By (\ref{3.14}) and (\ref{1.2}):
\begin{eqnarray}
\label{3.28} P(D_j^1)&=&(U;q)_j q^ju_1\cdot
P_{q^ju_1,q^{j+1}u_2,q^jU,q,i_0+i_1+i_2-j-1}\left(
\af{A}{i_0+i_1+i_2-j-1}{i_0-j}{i_1-1}{i_2}\right)\nonumber\\
&=&(U;q)_{i_0}u_1^{i_1}u_2^{i_2}
\qp{i_0+i_1+i_2-j-1}{i_0-j}{i_1-1}{i_2} q^{ji_1+(j+1)i_2}\,.
\end{eqnarray}
Analogously,
\begin{eqnarray}
\label{3.29} P(D_j^2)&=&(U;q)_j q^ju_1\cdot
P_{q^ju_1,q^ju_2,q^jU,q,i_0+i_1+i_2-j-1}\left(\af{A}
{i_0+i_1+i_2-j-1}{i_0-j}{i_1}{i_2-1}\right)\nonumber\\
&=&(U;q)_{i_0}u_1^{i_1}u_2^{i_2}
\qp{i_0+i_1+i_2-j-1}{i_0-j}{i_1}{i_2-1} q^{ji_1+ji_2}\,.
\end{eqnarray}
Substituting (\ref{3.28}), (\ref{3.29}), and (\ref{3.14}) into
(\ref{3.27}) and dividing by $(U;q)_{i_0}u_1^{i_1}u_2^{i_2}$, we
obtain
\begin{eqnarray}
\label{3.30} \qp{i_0+i_1+i_2}{i_0}{i_1}{i_2}&=&
\sum\limits_{j=0}^{i_0}
\qp{i_0+i_1+i_2-j-1}{i_0-j}{i_1-1}{i_2}q^{(i_1+i_2)j+i_2}+\nonumber\\
&+& \sum\limits_{j=0}^{i_0}
\qp{i_0+i_1+i_2-j-1}{i_0-j}{i_1}{i_2-1}q^{(i_1+i_2)j}\,.
\end{eqnarray}
If we replace $j$ in (\ref{3.30}) by $k=i_0-j$ and use identities
(\ref{3.26}), then we obtain (\ref{3.21}) after a simple
calculation.

3) The proof of (\ref{3.22}) is similar to those of (\ref{3.20})
and (\ref{3.21}). For $j,k\geq 0$, $j+k\leq i_0$, $a,b=1,2$ we
consider the events
\begin{eqnarray*}
E_{jk}^{ab}=\af{A}{i_0+i_1+i_2}{i_0}{i_1}{i_2}\cap
\{\!\!\!\!\!\!\!&&(\varepsilon_1,\ldots,\varepsilon_{i_0+i_1+i_2})
\in\Omega_{i_0+i_1+i_2}:
\varepsilon_1=\ldots=\varepsilon_{j}=0, \varepsilon_{j+1}=a,\\
&&\varepsilon_{i_0+i_1+i_2-k}=b,\varepsilon_{i_0+i_1+i_2-k+1}=\ldots
= \varepsilon_{i_0+i_1+i_2}=0\}\,.
\end{eqnarray*}
It can be  easily seen that the events $E_{jk}^{ab}$ are disjoint
with union $\af{A}{i_0+i_1+i_2}{i_0}{i_1}{i_2}$. Then we can write
\begin{eqnarray}
\label{3.31} P\left(\af{A}{i_0+i_1+i_2}{i_0}{i_1}{i_2}\right)=
\sum\limits_{j,k\geq 0\atop j+k\leq i_0}
\sum\limits_{a,b=1}^2P\left(E_{jk}^{ab}\right)\,.
\end{eqnarray}

Calculations similar to those that are made through proving of
(\ref{3.21}) show that
\be \label{3.32} P(E_{jk}^{11})
=(U;q)_{i_0}u_1^{i_1}u_2^{i_2}
\qp{i_0+i_1+i_2-j-k-2}{i_0-j-k}{i_1-2}{i_2}
q^{j(i_1+i_2)+i_0+i_2-j-k}\,, \ee \be \label{3.33} P(E_{jk}^{22})
=(U;q)_{i_0}u_1^{i_1}u_2^{i_2}
\qp{i_0+i_1+i_2-j-k-2}{i_0-j-k}{i_1}{i_2-2}
q^{j(i_1+i_2)+i_0+i_1-j-k}\,, \ee \be \label{3.34} P(E_{jk}^{21})
=(U;q)_{i_0}u_1^{i_1}u_2^{i_2}
\qp{i_0+i_1+i_2-j-k-2}{i_0-j-k}{i_1-1}{i_2-1}
q^{j(i_1+i_2)+i_0-j-k}\,, \ee \be \label{3.35} P(E_{jk}^{12})
=(U;q)_{i_0}u_1^{i_1}u_2^{i_2}
\qp{i_0+i_1+i_2-j-k-2}{i_0-j-k}{i_1-1}{i_2-1}
q^{j(i_1+i_2)+i_0++i_1+i_2-j-k-1}\,. \ee Substituting
(\ref{3.14}), (\ref{3.32}) -- (\ref{3.35}) into (\ref{3.31}) and
using (\ref{3.26}), we get after a division by
$(U;q)_{i_0}u_1^{i_1}u_2^{i_2}$ \be \label{3.36}
\begin{array}{rcl}
&&{\displaystyle\qp{i_0+i_1+i_2}{i_0}{i_1}{i_2}=
\left(1-q^{i_1+i_2-1}\right)\left(1-q^{i_1+i_2}\right)
\times}\nonumber\\
&&{\displaystyle\times\sum\limits_{j,k\geq 0\atop j+k\leq i_0}
\qp{i_0+i_1+i_2-j-k}{i_0-j-k}{i_1}{i_2} {q^{j(i_1+i_2)+i_0-j-k}
\over \left(1-q^{i_0+i_1+i_2-j-k-1}\right)
\left(1-q^{i_0+i_1+i_2-j-k}\right)}}\,.
\end{array}
\ee We obtain (\ref{3.22}) from (\ref{3.36}) after a change of
variables $j'=j$, $k'=i_0-j-k$ and a simple computation. $\Box$
\medskip

{\bf Remark 3.2.} It is clear that the arguments of Section~3 are
perfectly general. In order to generalize our arguments to the
case of $q$-polynomial coefficients
$\displaystyle\left[{i_0+\ldots+i_k}\atop
i_0,\ldots,i_k\right]_q$, $k\geq 3$, we can take a trial ${\cal
E}$ with $k+2$ outcomes $0, 1, 2,\ldots, k, \ast\,$. Analogously
to (\ref{3.1}), we assume that the probabilities of the outcomes
$0, 1, 2,\ldots, k$ in the first trial are equal to
$$
1-U\,,\,\, u_1\,,\,\, u_2\,,\,\,\ldots\,,\,\, u_k\,,
$$
respectively, where $u_i>0$ ($i=1, 2,\ldots, k$),
$U:=u_1+u_2+\ldots +u_k<1$. If $m$ is a positive integer, $i_0,
i_1,\ldots, i_k\geq 0$, $i_0+i_1+\ldots +i_k=m$, and if we know
that in the first $m$ trials the outcomes $0, 1, 2, \ldots, k$
have happened $i_0, i_1, i_2,\ldots, i_k$ times, respectively, we
assume that the probabilities of the outcomes $0, 1, 2,\ldots, k$
in the $(m+1)^{th}$ trial are equal to
$$
1-q^{i_0}U\,,\,\, q^{i_0}u_1\,,\,\,
q^{i_0+i_1}u_2\,,\,\,\ldots\,,\,\, q^{i_0+i_1+\ldots
+i_{k-1}}u_k\,,
$$ respectively.
\bigskip


\section{Stirling numbers of the second kind.}
\setcounter{equation}{0}
\medskip

Let $\Omega_n:=
\{\omega=(\varepsilon_1,\varepsilon_2,\ldots,\varepsilon_n):
\varepsilon_j=0\mbox{ or }1, j=1,2,\ldots,n\}$ be a set of all
sequences of the length $n$ with elements $0$ and $1$,
${\bf\beta}=\{\beta_j\}_{j=0}^{\infty}$ be a sequence of positive
numbers. We define a weight $w_n$ on $\Omega_n$ inductively. For
$n=1$, we set \be \label{4.1} w_1((0))=1,\qquad
w_1((1))=\beta_0\,. \ee Let $m>1$. We define the weight of a chain
of the length $m$ to be \be \label{4.2}
w_m((\varepsilon_1,\varepsilon_2,\ldots ,\varepsilon_{m-1},1))=
w_{m-1}((\varepsilon_1,\varepsilon_2,\ldots
,\varepsilon_{m-1}))\cdot\beta_j\,, \ee where $j=\#\{l: 1\leq
l\leq m-1,\varepsilon_l=0\}$, and \be \label{4.3}
w_m((\varepsilon_1,\varepsilon_2,\ldots ,\varepsilon_{m-1},0))=
w_{m-1}((\varepsilon_1,\varepsilon_2,\ldots
,\varepsilon_{m-1}))\,. \ee

In other words, the weight of a chain
$(\varepsilon_1,\varepsilon_2,\ldots
,\varepsilon_{m-1},\varepsilon_m)$ of the length $m$  equals  the
product of the weight of the chain
$(\varepsilon_1,\varepsilon_2,\ldots ,\varepsilon_{m-1})$ and that
of the element $\varepsilon_m,$  which is equal to $1$ if
$\varepsilon_m=0$ and to $\beta_j$ if $\varepsilon_m=1$ and
$\#\{k:1\leq k\leq m-1, \varepsilon_k=0\}=j$.

For every set $A\subset \Omega_n$, we define the weight $W_n(A)$
of $A$ to be \be \label{4.4} W_n(A):=\sum\limits_{\omega\in
A}w_n(\omega)\,. \ee It is evident from (\ref{4.4}) that the
additive property of the weight is valid: \be \label{4.5}
W_n(A\cup B)=W_n(A)+W_n(B)\,,\quad\mbox{ if } A\cap B=\emptyset\,.
\ee
For $n\geq 1$ and $0\leq k\leq n$ we denote \be \label{4.6}
\xi_{nk}:=\xi_{nk}({\bf\beta}):=W_n({\bf 0}(^n_k)) \,, \ee where
${\bf 0}(^n_k)=\{(\varepsilon_1,\varepsilon_2,\ldots
,\varepsilon_n)\in\Omega_n:\#\{l:1\leq l\leq n,
\varepsilon_l=0\}=k\}$ (see (\ref{2.4})). We denote also
$\xi_{00}(\beta)=1$, $\xi_{nk}(\beta)=0$ if $k<0$ or $k>n$. We see
that $\xi_{nk}$ is a polynomial in the variables $\beta_i$, $0\leq
i\leq k$, considering $\beta_j$ as independent variables.
\medskip

{\bf Definition 4.1.} Polynomials $\xi_{nk}$ are said to be {\it
Stirling polynomials of the second kind} generated by the sequence
$\beta$.
\medskip

The following theorem gives a recurrence relation for polynomials
$\xi_{nk}$.
\medskip

{\bf Theorem 4.1.} {\it If $n\in{\bf N}$ and $0\leq k\leq n$, then
\be \label{4.7} \xi_{nk}=\xi_{n-1,k-1}+\xi_{n-1,k}\beta_k\,. \ee }

{\bf Proof.} For $j=0, 1$, we denote $A^j:={\bf
0}(^n_k)\cap\{(\varepsilon_1,\ldots, \varepsilon_n)\in
\Omega_n:\varepsilon_n=j\}$. From (\ref{4.2}) and (\ref{4.3}) it
follows that
\begin{eqnarray*}
\xi_{nk}&=&W_n({\bf 0}(^n_k))=W_n(A^0)+W_n(A^1)\\
&=&W_{n-1}({\bf 0}(^{n-1}_{k-1}))\cdot 1+
W_{n-1}({\bf 0}(^{n-1}_k))\cdot\beta_k\\
&=&\xi_{n-1,k-1}+\xi_{n-1,k}\beta_k\,,
\end{eqnarray*}
This proves (\ref{4.7}). $\Box$
\medskip

For every positive integer $l$ and a sequence ${\bf\beta}:=
\{\beta_j\}_{j=0}^{\infty}$, we denote ${\bf\beta}^{(l)}:=
\{\beta_{l+j}\}_{j=0}^{\infty}$. The $W_n^{(l)}$  will denote the
weight on $\Omega_n$ generated by the sequence ${\bf\beta}^{(l)}$.
\medskip

{\bf Definition 4.2.} Polynomials
$$
\xi_{nk}^{(l)}:=\xi_{nk}({\bf\beta}^{(l)})\,,\quad
(n=1,2,\ldots\,,\,\,\,\, 0\leq k\leq n)\,;\quad
\xi_{00}^{(l)}:=1\,,
$$
in the variables $\beta_l, \beta_{l+1},\ldots$ are said to be {\it
associated of the rank $l$ } with polynomials
$\xi_{nk}({\bf\beta})$.
\medskip

The following theorem gives a relation that includes $\xi_{nk}$
and $\xi_{nk}^{(1)}$.
\medskip

{\bf Theorem 4.2.} {\it For all $n\geq 1$ and $0\leq k\leq n$ the
following recurrence relation holds: \be \label{4.8}
\xi_{nk}=\xi_{n-1,k-1}^{(1)}+\xi_{n-1,k}\beta_0\,. \ee }
\medskip

{\bf Proof.} For $j=0, 1$ we denote $B^j:={\bf
0}(^n_k)\cap\{(\varepsilon_1,\ldots, \varepsilon_n)\in
\Omega_n:\varepsilon_1=j\}$. Using (\ref{4.2}) and (\ref{4.3}), we
get
\begin{eqnarray*}
\xi_{nk}&=&W_n({\bf 0}(^n_k))= W_n(B^0)+
W_n(B^1)\\
&=&1\cdot W_{n-1}^{(1)}({\bf 0}(^{n-1}_{k-1}))+
\beta_0\cdot W_{n-1}({\bf 0}(^{n-1}_k))\\
&=&\xi_{n-1,k-1}^{(1)}+\xi_{n-1,k}\beta_0\,.\quad\Box
\end{eqnarray*}

Let us consider a particular case. We put \be \label{4.9}
\beta_j:=j \quad\mbox{ for all }\quad j\geq 0\,. \ee (Therefore,
$w_n(\omega)=0$ for every chain
$\omega=(\,\varepsilon_1,\varepsilon_2,\ldots,\varepsilon_n\,)$
such that $\varepsilon_1=1$.) Then we get numbers
${\tilde\xi}_{nk}:=\xi_{nk}(\{j\}_{j=0}^{\infty})$ satisfying the
following recurrence relation (see (\ref{4.7})) \be \label{4.10}
\tilde\xi_{nk}=\tilde\xi_{n-1,k-1}+\tilde\xi_{n-1,k}\cdot k \ee
and conditions $\tilde\xi_{00}=1$, $\tilde\xi_{nk}=0$ if $k<0$ or
$k>n$. The theorem below follows directly from the definition of
Stirling numbers of the second kind $\displaystyle\stt{n}{k}$ (see
(\ref{1.3})).
\medskip

{\bf Theorem 4.3.} {\it Let $n\in{\bf N}$, $0\leq k\leq n$. Then
\be \label{4.11} \stt{n}{k}=W_n({\bf 0}(^n_k))\,, \ee where $W_n$
denotes the weight on $\Omega_n$ generated by the sequence
$\beta=\{j\}_{j=0}^{\infty}$ with the help of $(\ref{4.1})$
--- $(\ref{4.4})$. }
\medskip

In the following theorem we give new proofs of some known facts
(see, for example, [3], formulas (6.20), (6.22)). Our proofs are
based on Theorem~4.3.
\medskip

{\bf Theorem 4.4.} 1) {\it If $n\in{\bf N}$ and $0\leq m\leq n$,
then \be \label{4.12}
\stt{n}{m}=\sum\limits_{l=m}^n\stt{l-1}{m-1}m^{n-l}\,. \ee } 2)
{\it If $n,m\in{\bf N}_0$, then \be \label{4.13}
\stt{n+m+1}{m}=\sum\limits_{k=0}^mk\stt{n+k}{k}\,. \ee }

{\bf Proof.} 1) We denote
$$
F_l:={\bf 0}(^n_m)\cap\{(\varepsilon_1,\ldots , \varepsilon_n):
\varepsilon_l=0, \varepsilon_{l+1}=\varepsilon_{l+2}= \ldots
=\varepsilon_n=1\}
$$
for every $l=m, m+1,\ldots ,n$. It is evident that these sets form
a partition of ${\bf 0}(^n_m)$. We calculate $W_n(F_l)$. Every
chain $\omega$ from the set $F_l$ has the form
$\omega=(\omega^{\prime}, 0,1_{(n-l)} )$, where
$\omega^{\prime}\in{\bf 0}(^{l-1}_{m-1})$ and $j_{(k)}$ denotes
the sequence $\displaystyle\underbrace{j,j,\ldots,j}_{k}$.
Therefore, by (\ref{4.2}), (\ref{4.3}), (\ref{4.11}), and
(\ref{4.9}) \be \label{4.14} W_n(F_l)=W_{l-1}({\bf
0}(^{l-1}_{m-1}))\cdot 1\cdot \beta_m^{n-l}= \stt{l-1}{m-1}\cdot
m^{n-l}\,. \ee Inserting (\ref{4.11}) and (\ref{4.14}) into
 $W_n({\bf 0}(^n_m))=\sum_{l=m}^nW_n(F_l)$, we obtain
(\ref{4.12}).

2) For every $k=0,1,\ldots ,m$, we introduce the sets
$$
H_k:= {\bf 0}(^{n+m+1}_m)\cap \{(\varepsilon_1,\ldots ,
\varepsilon_{n+m+1}): \varepsilon_{n+k+1}=1, \varepsilon_{n+k+2}=
\ldots =\varepsilon_{n+m+1}=0\}\,.
$$
It is evident that these sets form a partition of the set ${\bf
0}(^{n+m+1}_m)$. We now evaluate $W_{n+m+1}(H_k)$. Every chain
$\omega\in H_k$ has the form $\omega= (\omega^{\prime},1,
0_{m-k})$, where $\omega^{\prime}\in{\bf 0}(^{n+k}_k)$. Therefore
\be \label{4.15} W_{n+m+1}(H_k)=W_{n+k}({\bf
0}(^{n+k}_k))\cdot\beta_k =\stt{n+k}{k}\cdot k\,. \ee Inserting
(\ref{4.11}) and (\ref{4.15}) into
$$ W_{n+m+1}({\bf
0}(^{n+m+1}_m))=\sum_{k=0}^mW_{n+m+1}(H_k)\,,
$$
we get (\ref{4.13}). $\Box$
\medskip

The following theorem gives relations between Stirling numbers of
the second kind and associated ones with them.
\medskip

{\bf Theorem 4.5.} 1) {\it If $n, m\in{\bf N}$ and $1\leq m\leq
n$, then \be \label{4.16}
\stt{n}{m}=\sum\limits_{j=1}^mj\stt{n-j-1}{m-j}^{(j)}\,. \ee } 2)
{\it If $n, \nu\in{\bf N}$, $1\leq\nu\leq n-1$, $0\leq m\leq n$,
then \be \label{4.17}
\stt{n}{m}=\sum\limits_{k=0}^{\nu}\stt{\nu}{k}\stt{n-\nu}{m-k}^{(k)}\,.
\ee }
\medskip

{\bf Proof.} 1) For $j=0, 1, 2,\ldots ,m$, we consider the sets
$$
G_j:={\bf 0}(^n_m)\cap\{(\varepsilon_1,\ldots ,
\varepsilon_n)\in\Omega_n: \varepsilon_1=\ldots
=\varepsilon_j=0,\varepsilon_{j+1}=1\}.
$$
Evidently, $G_j$ form a partition of ${\bf 0}(^n_m)$. We evaluate
$W_n(G_j)$. Every chain $\omega\in G_j$ has the form $\omega=
(0_{(j)},1,\omega^{\prime})$, where $\omega^{\prime}\in{\bf
0}(^{n-j-1}_{m-j})$. Therefore, by the definition of the weight
$W_n$, we have \be \label{4.18} W_n(G_j)= \beta_j\cdot
W_{n-j-1}^{(j)}({\bf 0}(^{n-j-1}_{m-j}))=
j\cdot\stt{n-j-1}{m-j}^{(j)}\,. \ee We recall that $W^{(j)}$ is a
weight generated by the sequence $\{j+i\}_{i=0}^{\infty}$ and
$\displaystyle\stt{a}{b}^{(j)}:=W_n^{(j)}({\bf 0}(^a_b))$ are
numbers associated with Stirling ones of the rank $j$. Inserting
(\ref{4.11}) and (\ref{4.18}) into  $W_n({\bf
0}(^n_m))=\sum_{j=0}^mW_n(G_j)$, we get (\ref{4.16}).

2) For $k=0,1,2,\ldots,\nu$, we consider the sets $R_k:={\bf
0}(^n_m)\cap{\bf 0}(^{\nu}_k)$. Evidently, $R_k$ form a partition
of ${\bf 0}(^n_m)$. We evaluate $W_n(R_k)$. Every chain $\omega\in
R_k$ has the form
$\omega=(\omega^{\prime},\omega^{\prime\prime})$, where
$\omega^{\prime}\in{\bf 0}(^{\nu}_k)$, $\omega^{\prime\prime}\in
{\bf 0}(^{n-\nu}_{m-k})$. By the definition of the weight $W_n$,
we have \be \label{4.19} W_n(R_k)= W_{\nu}({\bf
0}(^{\nu}_k))W_{n-\nu}^{(k)}({\bf 0}(^{n-\nu}_{m-k}))=
\stt{\nu}{k}\stt{n-\nu}{m-k}^{(k)}\,. \ee Inserting (\ref{4.11})
and (\ref{4.19}) into  $W_n({\bf
0}(^n_m))=\sum_{k=0}^{\nu}W_n(R_k)$, we obtain (\ref{4.17}).
$\Box$
\medskip

{\bf Remark 4.1.} It is easy to generalize Theorems~4.4 and 4.5 to
the case of an arbitrary sequence $\{\beta_j\}$.
\bigskip


\section{Stirling numbers of the first kind.}
\setcounter{equation}{0}
\bigskip

As in Section~4, let
$\Omega_n=\{(\varepsilon_1,\varepsilon_2,\ldots
,\varepsilon_n):\varepsilon_j=0, 1; 1\leq j\leq n\}$,
${\bf\gamma}=\{\gamma_j\}_{j=1}^{\infty}$ be a sequence of
positive numbers. We introduce a weight $w_n$ on the set
$\Omega_n$ inductively. We put for $n=1$, \be \label{5.1}
w_1((0))=1,\qquad w_1((1))=\gamma_1\,, \ee and for $m>1$, \be
\label{5.2}
\begin{array}{rcl}
w_m((\varepsilon_1,\varepsilon_2,\ldots ,\varepsilon_{m-1},0))&=&
w_{m-1}((\varepsilon_1,\varepsilon_2,\ldots
,\varepsilon_{m-1}))\,,\\
w_m((\varepsilon_1,\varepsilon_2,\ldots ,\varepsilon_{m-1},1))&=&
w_{m-1}((\varepsilon_1,\varepsilon_2,\ldots,
\varepsilon_{m-1}))\cdot\gamma_m\,.
\end{array}
\ee We define $W_n(A)$ for all $A\subset\Omega_n$ as in
(\ref{4.4}).

For every sequence $\gamma=\{\gamma_j\}_{j=1}^{\infty}$, for all
$n\in{\bf N}$ and integer $k$, $0\leq k\leq n$,  we define
polynomials $\eta_{nk}$ in variables $\gamma_j$ as follows
\be\label{5.14} \eta_{nk}:=\eta_{nk}(\gamma):= W_n({\bf
0}(^n_k))\,. \ee For every sequence $\gamma$ we define also
$\eta_{00}(\gamma)=1$, $\eta_{nk}(\gamma)=0$ if $k<0$ or $k>n$.
For each nonnegative integer $n$ and $0\leq k\leq n$, $\eta_{nk}$
is a polynomial in the variables $\gamma_i$, $1\leq i\leq n$.
\medskip

{\bf Definition 5.1.} We say that polynomials $\eta_{nk}$ are {\it
Stirling polynomials of the first kind} generated by the sequence
$\gamma$.
\medskip

Polynomials $\eta_{nk}(\gamma)$ satisfy the following recurrence
relation.
\medskip

{\bf Theorem 5.1.} {\it Let $n\in{\bf N}$, $0\leq k\leq n$, and
$\eta_{nk}(\gamma)$ be polynomials in $\gamma_1,\ldots ,\gamma_n$
defined by $(\ref{5.14})$. Then \be \label{5.3} \eta_{nk}(\gamma)=
\eta_{n-1,k-1}(\gamma)+\eta_{n-1,k}(\gamma)\cdot\gamma_n\,. \ee }

This is an analogue of Theorem~4.1 and the proof is the same.
$\Box$
\medskip

As in the previous section, we write
$\gamma^{(l)}:=\{\gamma_{l+j}\}_{j=1}^{\infty}$ for every
$l\in{\bf N}$.
\medskip

{\bf Definition 5.2.} Polynomials
$\eta_{nk}^{(l)}:=\eta_{nk}(\gamma^{(l)})$ in the variables
$\gamma_{l+1},\gamma_{l+2},\ldots$ are said to be {\it associated
ones of the rank $l$} with polynomials $\eta_{nk}(\gamma)$.
\medskip

We consider now  a particular case: \be \label{5.4} \gamma_j:=j-1
\quad\mbox{ for all }\quad j\geq 1\,. \ee (Therefore,
$w_n(\omega)=0$ for every chain
$\omega=(\,\varepsilon_1,\varepsilon_2,\ldots,\varepsilon_n\,)$
such that $\varepsilon_1=1$.) We get a set of numbers
$\tilde\eta_{nk}:=\eta_{nk}(\{j-1\}_{j=1}^{\infty})$, which
satisfy the recurrence relation
$$
\tilde\eta_{nk}= \tilde\eta_{n-1,k-1}+\tilde\eta_{n-1,k}\cdot
(n-1)
$$
and conditions $\tilde\eta_{00}=1$, $\tilde\eta_{nk}=0$ if $k<0$
or $k>n$. These numbers are called as {\it Stirling numbers of the
first kind} and are denoted by $\displaystyle\st{n}{k}$ (see
(\ref{1.3})).

As a result, we can derive the following statement.
\medskip

{\bf Theorem 5.2.} {\it Let $n\in{\bf N}$, $0\leq k\leq n$. Then
\be \label{5.5} \st{n}{k}=W_n({\bf 0}(^n_k))\,, \ee where $W_n$ is
the weight on $\Omega_n$ generated by the sequence
$\gamma=\{j-1\}_{j=1}^{\infty}$ with the help of $(\ref{5.1})$,
$(\ref{5.2})$, $(\ref{4.4})$ } \medskip

Using (\ref{5.5}), we will give very simple proofs of two known
facts (see, for example, [3], formulas (6.21), (6.23)).
\bigskip

{\bf Theorem 5.3.} 1) {\it If $n\in{\bf N}$ and $m$ is an integer
such that $0\leq m\leq n$, then \be \label{5.6}
\st{n}{m}=\sum\limits_{l=m}^n\st{l-1}{m-1}l(l+1)(l+2)\ldots(n-1)\,.
\ee } 2) {\it If $n, m\in{\bf N}$, then \be \label{5.7}
\st{n+m+1}{m}=\sum\limits_{k=0}^m(n+k)\st{n+k}{k}\,. \ee }

{\bf Proof.} The proof is very similar to that of Theorem~4.4. To
prove (\ref{5.6}) and (\ref{5.7}), we consider the sets $F_l$,
$l=m,m+1,\ldots,n,$ and $H_k$, $k=0,1,\ldots,m,$ introduced in the
proof of the first and the second parts of Theorem~4.4,
respectively. We find from (\ref{5.1}) and (\ref{5.2}) that
\begin{eqnarray*}
W_n(F_l)&=&W_{l-1}({\bf 0}(^{l-1}_{m-1}))
\gamma_{l+1}\gamma_{l+2}\cdot\ldots\cdot\gamma_n\nonumber\\
&=&\st{l-1}{m-1}l(l+1)(l+2)\ldots(n-1)\,,\\
W_{n+m+1}(H_k)&=&W_{n+k}({\bf 0}(^{n+k}_k))\cdot\gamma_{n+k+1} =
\st{n+k}{k}\cdot (n+k)\,.
\end{eqnarray*}
Repeating the reasoning from the proof of Theorem~4.4, we obtain
(\ref{5.6}), (\ref{5.7}). $\Box$
\medskip

The following theorem is an analogue of Theorem~4.5.
\medskip

{\bf Theorem 5.4.} 1) {\it If $n, m\in{\bf N}$, $1\leq m\leq n$,
then \be \label{5.10}
\st{n}{m}=\sum\limits_{j=1}^mj\st{n-j-1}{m-j}^{(j)}\,. \ee } 2)
{\it If $n, \nu, m\in{\bf N}$, $1\leq \nu,m\leq n$, then \be
\label{5.11}
\st{n}{m}=\sum\limits_{k=0}^{\nu}\st{\nu}{k}\st{n-\nu}{m-k}^{(\nu)}\,.
\ee }
\medskip

{\bf Proof.} The proof  is similar to that of Theorem~4.5. We
consider the sets $G_j$ ($1,2,\ldots,m$), $R_k$
($k=0,1,\ldots,\nu$)introduced there. In our case the weights of
these sets are equal to
\begin{eqnarray*}
&&W_{n+m+1}(G_j)=\gamma_{j+1}\cdot
\st{n-j-1}{m-j}^{(j)}= j\cdot\st{n-j-1}{m-j}^{(j)}\,,\\
&&W_{n+m+1}(R_k)= W_{\nu}({\bf 0}(^{\nu}_k))W_{n-\nu}^{(\nu)}({\bf
0}(^{n-\nu}_{m-k}))= \st{\nu}{k}\st{n-\nu}{m-k}^{(\nu)}\,.
\end{eqnarray*}
The theorem is now immediate. $\Box$
\medskip

{\bf Remark 5.1.} It is easy to generalize Theorems~5.3 and 5.4 to
the case of an arbitrary sequence $\{\gamma_j\}$.
\bigskip


\section{Euler numbers.}
\setcounter{equation}{0}
\bigskip

As before, let $\Omega_n$ be the set of all sequences of the
length $n$ with elements $0$ and $1$,
$\{\alpha_j\}_{j=1}^{\infty}$ and $\{\beta_j\}_{j=0}^{\infty}$ be
two sequences of positive numbers. Let us introduce a weight on
$\Omega_n$ by induction on $n$. For $n=1$, we set \be \label{6.1}
w_1((0))=\alpha_1\,,\qquad w_1((1))=\beta_0\,. \ee Let $m>1$. We
define the weight of a chain of the length $m$ as follows \be
\label{6.2}
\begin{array}{rl}
&w_m((\varepsilon_1,\varepsilon_2,\ldots ,\varepsilon_{m-1},0))=
w_{m-1}((\varepsilon_1,\varepsilon_2,\ldots ,\varepsilon_{m-1}))
\cdot\alpha_{m-k}\,,
\\
&w_m((\varepsilon_1,\varepsilon_2,\ldots ,\varepsilon_{m-1},1))=
w_{m-1}((\varepsilon_1,\varepsilon_2,\ldots,
\varepsilon_{m-1}))\cdot\beta_k\,,
\end{array}
\ee where $k=\#\{j:1\leq j\leq m-1, \varepsilon_j=0\}$. In other
words, the weight of a chain of the length $m$ equals the product
of the weight of the chain consisting of the first $m-1$ terms of
a given one and of the weight of the $m^{th}$ term which is equal
to $\alpha_{m-k}$, if this term is $0$, and $\beta_k$, if it is
$1$ and if $k$ terms are $0$ among the first $m-1$ ones of a given
chain. Evidently, definition (\ref{6.1}) is consistent with
definition (\ref{6.2}), that is (\ref{6.1}) follows from
(\ref{6.2}) if we take $m=1$ and if we assume that the first term
at the right-hand side of both equalities (\ref{6.2}) equals $1$.
As before, we define the weight $W_n(A)$ of a set $A\subset
\Omega_n$ by the formula (\ref{4.4}).

For every $n\in{\bf N}$ and integer $k$ such that $0\leq k\leq n$,
we define \be\label{6.15} \zeta_{nk}:=\zeta_{nk}(\alpha ,\beta):=
W_n({\bf 0}(^n_k))\,. \ee By definition we put $\zeta_{00}(\alpha,
\beta)=1$ and $\zeta_{nk}(\alpha, \beta)=0$ whenever $k<0$ or
$k>n$. It is evident that $\zeta_{nk}$ are polynomials in the
variables $\alpha_i$, $\beta_j$ (if we consider $\alpha_i$,
$\beta_j$ as independent variables).
\medskip

{\bf Definition 6.1.} Polynomials $\zeta_{nk}(\alpha,\beta)$ are
said to be {\it Euler polynomials}, generated by sequences
$\alpha$ and $\beta$.
\medskip

The following theorem gives a recurrence relation for the
polynomials $\zeta_{nk}$.
\medskip

{\bf Theorem 6.1.} {\it Let $n\in{\bf N}$, $0\leq k\leq n$, and
$\zeta_{nk}(\alpha, \beta)$ be polynomials defined by
$(\ref{6.15})$. Then \be \label{6.3} \zeta_{nk}(\alpha,
\beta)=\zeta_{n-1,k-1}(\alpha,
\beta)\alpha_{n-k+1}+\zeta_{n-1,k}(\alpha, \beta)\beta_k\,. \ee }
\medskip

{\bf Proof.} The proof is  analogous to that of Theorems~4.1 and
5.1. We only note that if $A_0$ and $A_1$ are defined as in the
proof of Theorem~4.1, then
$$
W_n(A_0)=W_{n-1}({\bf
0}(^{n-1}_{k-1}))\cdot\alpha_{n-(k-1)}\,,\qquad
W_n(A_1)=W_{n-1}({\bf 0}(^{n-1}_k))\cdot\beta_k\,.\Box
$$
\medskip

{\bf Definition 6.2.} For every $\nu\in{\bf N}$ and integer $\mu$
such that $0\leq\mu\leq\nu$, we define the polynomials
$$
\zeta_{nk}^{(\nu,\mu)}(\alpha,\beta):=
\zeta_{nk}(\alpha^{(\nu-\mu)},\beta^{(\mu)})
$$
and call them {\it polynomials associated with the polynomials
$\zeta_{nk}(\alpha,\beta)$  of rank $(\nu,\mu)$}.   They are
polynomials in variables $\alpha_{\nu-\mu+1},
\alpha_{\nu-\mu+2},\ldots$, $\beta_{\mu}, \beta_{\mu+1},\ldots$
(We recall that if $\delta=\{\delta_j\}_{j=j_0}^{\infty}$ is a
sequence and $l$ is a nonnegative integer, then we denote
$\delta^{(l)}:=\{\delta_{l+j}\}_{j=j_0}^{\infty}$.)
\medskip

We consider a particular case. Let $$ \alpha_l=l-1 \mbox{ for all
} l\geq 1\,,\qquad \beta_k=k+1 \mbox{ for all } k\geq 0\,. $$ We
obtain a set of numbers $\tilde\zeta_{nk}:=
\zeta_{nk}(\{l-1\}_{l=1}^{\infty},\{k+1\}_{k=0}^{\infty})$,
($n\in{\bf N}_0$, $0\leq k\leq n$) such that \be\label{6.4}
\tilde\zeta_{nk}=\tilde\zeta_{n-1,k-1}\cdot
(n-k)+\tilde\zeta_{n-1,k}\cdot (k+1)\,, \ee and
$\tilde\zeta_{00}=1$, $\tilde\zeta_{nk}=0$ whenever $k<0$ or
$k>n$. These numbers are called {\it Euler numbers} and are
denoted by $\displaystyle\eu{n}{k}$ (see (\ref{1.3})). By
(\ref{1.3}) and (\ref{6.4}), the following theorem holds.
\medskip

{\bf Theorem 6.2.} {\it Let $n\in{\bf N}$ and $0\leq k\leq n$.
Then \be \label{6.5} \eu{n}{k}=W_n({\bf 0}(^n_k))\,, \ee where the
weight $W_n$ on $\Omega_n$ is generated by the sequences
$\alpha=\{l-1\}_{l=1}^{\infty}$, $\beta=\{k+1\}_{k=0}^{\infty}$ by
means of $(\ref{6.1})$, $(\ref{6.2})$, and $(\ref{4.4})$. }
\medskip

The following theorem is an analogue of Theorems~4.4 and 5.3.
\medskip

{\bf Theorem 6.3.} 1) {\it If $n, m\in{\bf N}_0$, $0\leq m\leq n$,
then \be \label{6.6}
\eu{n}{m}=\sum\limits_{l={m+1}}^n\eu{l-1}{m-1}(l-m)(m+1)^{n-l}\,.
\ee $2)$ If $n, m\in{\bf N}_0$, then \be \label{6.7}
\eu{n+m+1}{m}=\sum\limits_{k=0}^m\eu{n+k}{k}(k+1)(n+1)^{m-k}\,.
\ee }

{\bf Proof.} For all $l=m, m+1,\ldots ,n$ and $k=0,1,\ldots,m$ we
introduce  the sets $F_l$ and $H_k$ in the same way as in the
proof of Theorem~4.4. By Theorem~6.2 we have
\begin{eqnarray*}
&&W_n(F_l)=W_{l-1}({\bf 0}(^{l-1}_{m-1}))\alpha_{l-(m-1)}
\beta_m^{n-l}=
\eu{l-1}{m-1}(l-m)(m+1)^{n-l}\,,\\
&&W_{n+m+1}(H_k)=W_{n+k}({\bf 0}(^{n+k}_k))\beta_k
\alpha_{n+2}^{m-k} =\eu{n+k}{k}(k+1)(n+1)^{m-k}\,.
\end{eqnarray*}
The theorem is now immediate. $\Box$
\medskip

{\bf Remark 6.1.} {\it For all $n\in{\bf N}$ and $m$ such that
$0\leq m\leq n-1$, the following identity is valid: \be
\label{6.10} \eu{n}{m}=\eu{n-1}{m}^{(1,0)}. \ee }
\medskip

For every $j=0,1, 2,\ldots ,m$, we set
$$
Q_j:={\bf 0}(^n_m)\cap\{\omega=(0_{(j)},1,\omega^{\prime}):
\omega^{\prime}\in{\bf 0}(^{n-j-1}_{m-j})\}\,.
$$
Therefore \be \label{6.11}
\begin{array}{rcl}
W_n(Q_j)&=&\underbrace{\alpha_{1-0}\cdot\alpha_{2-1}\cdot
\ldots\cdot\alpha_{j-(j-1)}}_{j}\cdot\beta_j\cdot
W_{n-j-1}^{(j+1,j)}({\bf 0}(^{n-j-1}_{m-j}))\\
&=&{\displaystyle\alpha_1^j\beta_j\cdot\eu{n-j-1}{m-j}^{(j+1,j)}=
\cases{0,& if $j\geq 1$;\cr \eu{n-1}{m}^{(1,0)},& if $j=0$\,.\cr}}
\end{array}
\ee Inserting (\ref{6.5}) and (\ref{6.11}) into  $W_n({\bf
0}(^n_m))= \sum_{j=0}^mW_n(Q_j)$, we obtain (\ref{6.10}).
\medskip

{\bf Theorem 6.4.} 1) {\it If $n,m\in{\bf N}$, $1\leq m\leq n$,
then \be \label{6.12}
\eu{n}{m}=\sum\limits_{j=1}^{n-m}j\eu{n-j-1}{m-j}^{(j+1,1)}\,. \ee
}

2) {\it If $n, m, \nu\in{\bf N}_0$, $0\leq m\leq n$, $1\leq\nu\leq
n-1$, then \be \label{6.13}
\eu{n}{m}=\sum\limits_{k=0}^{\nu}\eu{\nu}{k}\eu{n-\nu}{m-k}^{(\nu,k)}\,.
\ee }

{\bf Proof.} 1) Just as in the proof of Theorem~4.5, we consider
the sets $G_j$ ($j=0,1, \ldots ,m$). We have in our case \be
\label{6.14} W_n(G_j)=\beta_0^j\alpha_{j+1}
W_{n-j-1}^{(j+1,1)}({\bf 0}(^{n-j-1}_{m-1}))=
j\cdot\eu{n-j-1}{m-1}^{(j+1,1)}\,. \ee

2) For $k= 0,1,2,\ldots,\nu$, we set $R_k:={\bf 0}(^n_m)\cap{\bf
0}(^{\nu}_k)$. Every element $\omega$ of $R_k$ has the form
$\omega=(\omega^{\prime},\omega^{\prime\prime})$, where
$\omega^{\prime}\in{\bf 0}(^{\nu}_k)$,
$\omega^{\prime\prime}\in{\bf 0}(^{n-\nu}_{m-k})$. Therefore
$$
W_n({\bf 0}(^n_m))=\sum\limits_{k=0}^{\nu}W_n(R_k)=
\sum\limits_{k=0}^{\nu}W_{\nu}({\bf 0}(^{\nu}_k))
W_{n-\nu}^{(\nu,k)}({\bf 0}(^{n-\nu}_{m-k}))=
\sum\limits_{k=0}^{\nu}\eu{\nu}{k} \eu{n-\nu}{m-k}^{(\nu,k)}\,.
$$
\medskip

{\bf Remark 6.2.} It is easy to generalize Theorems~6.3 and 6.4 to
the case of arbitrary sequences $\{\alpha_j\}$ and $\{\beta_j\}$.


\bigskip

\begin{thebibliography}{99}

\bibitem{1} G. Andrews,
{\it The Theory of Partitions}, Addison-Wesley Publishing Company,
London, Amsterdam, 1976.
\bibitem{2} W. Feller,
{\it An Introduction to Probability Theory and its Applications},
Vol.~1, Third Edition, John Wiley \& Sons, New York, 1970.
\bibitem{3} R.L.Graham, D.E.Knuth, O.Patashnik,
{\it Concrete Mathematics}, Addison-Wesley Publishing Company,
London, Amsterdam, 1998.
\bibitem{4} A. Il'inskii, S. Ostrovska,
{\it Convergence of generalized Bernstein polynomials}, Journal of
Approximation Theory, {\bf 116} (2002), 100 -- 112.
\bibitem{5} R.P.Stanley,
{\it Enumerative Combinatorics}, Vol. 1, Wadsworth \& Brooks,
1986.
\end{thebibliography}
\end{document}